\documentclass[11pt]{amsart}

\usepackage[margin=1.08in]{geometry}
\usepackage{amsmath,amssymb,amsthm,mathtools}
\usepackage[expansion=false]{microtype}
\usepackage{hyperref}
\hypersetup{hidelinks}
\usepackage{enumitem}
\usepackage{tikz}

\numberwithin{equation}{section}

\newtheorem{theorem}{Theorem}[section]
\newtheorem{proposition}[theorem]{Proposition}
\newtheorem{lemma}[theorem]{Lemma}
\newtheorem{corollary}[theorem]{Corollary}
\newtheorem{conjecture}[theorem]{Conjecture}

\newcommand{\HH}{\mathfrak H}
\newcommand{\QQ}{\mathbb Q}
\newcommand{\ZZ}{\mathbb Z}
\newcommand{\CC}{\mathbb C}
\newcommand{\Qbar}{\overline{\mathbb Q}}
\newcommand{\PP}{\mathbb P}
\newcommand{\trdeg}{\operatorname{trdeg}}

\newcommand{\RR}{\mathbb R}
\newcommand{\aclC}[1]{\overline{#1}^{\,\CC}}

\begin{document}
\title{Zeros of $E_2$, quasi-periods, and Schwarz maps}
\author{Abdellah Sebbar}
\address{Department of Mathematics, University of Ottawa}
\email{asebbar@uottawa.ca}


\begin{abstract}
We study the zeros of the Eisenstein series $E_2$, their distribution
near rational cusps, and the values of modular functions at these zeros.
Using Watson's inversion formula, we prove that the number of zeros in
$0\leq\operatorname{Re}\tau<1$ with imaginary part at least $Y$ satisfies
\[
 N(Y)=\frac{1}{2\pi Y}
      +O\!\left(Y^{-1/2}\log\frac1Y\right)
 \qquad(Y\to0^+).
\]
We also determine the limiting distribution of the scaled displacements
from the cusps and prove that two zeros whose difference is rational
differ by an integer. At every zero $\rho$, classical transcendence
theorems imply that both $j(\rho)$ and $\lambda(\rho)$ are
transcendental. We use this fact to obtain an exact multiplicity formula
for arithmetic quasimodular forms at arbitrary congruence level and to
show that all ordinary derivatives of arithmetic modular functions on
congruence subgroups, evaluated at $\rho$, are algebraic over
$\overline{\mathbb Q}(j(\rho),j'(\rho))$. We then relate Watson's
formula to the Schwarz map of a rational transform of the Legendre
hypergeometric equation. Finally, we give conditional results on the
transcendence of $\rho$ under modular Schanuel and period conjectures,
and determine why the associated analytic subgroup fails the rank
condition in Waldschmidt's theorem.
\end{abstract}

\maketitle
\setcounter{tocdepth}{1}
\begingroup
\footnotesize
\tableofcontents
\endgroup

\section{Introduction}

For the normalized lattice $\Lambda_\tau=\ZZ+\ZZ\tau$, the first
quasi-period is
\[
 \eta_1=\frac{\pi^2}{3}E_2(\tau),\qquad
 E_2(\tau)=1-24\sum_{n\ge1}\sigma_1(n)e^{2\pi i n\tau}.
\]
Thus $E_2(\tau)=0$ if and only if $\eta_1=0$. Legendre's relation gives
\begin{equation}\label{eq:intro-h}
 h(\tau)=\frac{\eta_2}{\eta_1}
 =\tau+\frac{6}{\pi iE_2(\tau)}.
\end{equation}
The poles of $h$ are the zeros of $E_2$. Watson
\cite[\S\S4--6, pp.~137--149]{Watson1952} parametrized these zeros
by rational cusps and obtained a convergent inversion formula.
Imamo\u{g}lu, Jermann and T\'oth
\cite[Theorems~1.1--1.2]{ImamogluJermannToth2014} proved that each
Ford circle contains exactly one simple zero. Wood and Young
\cite[Theorems~3.1, 5.1 and 6.2]{WoodYoung2014} independently
obtained localization estimates and a bijection between rational
values of $h$ and the zeros.

The zeros in one translation strip are indexed by reduced cusps and have the form
\begin{equation}\label{eq:intro-Watson}
 \rho_{m/n}=\frac mn+\frac1{n^2}\mathcal G\!\left(\frac{m'}n\right),
 \qquad mm'\equiv-1\pmod n,
\end{equation}
where $\mathcal G$ is one-periodic and real analytic.

\begin{theorem}[Counting and distribution]\label{thm:intro-main}
Let $N(Y)$ count the zeros $\rho$ with
$0\leq\operatorname{Re}\rho<1$ and $\operatorname{Im}\rho\geq Y$.
Then, as $Y\to0^+$,
\[
 \int_0^1\mathcal G(x)\,dx=\frac{\pi i}{6},\qquad
 N(Y)=\frac1{2\pi Y}+O\!\left(Y^{-1/2}\log\frac1Y\right).
\]
The probability measures
\[
 \frac1{\varphi(n)}
 \sum_{\substack{0\le m<n\\(m,n)=1}}
 \delta_{\,n^2(\rho_{m/n}-m/n)}
\]
converge weakly to $\mathcal G_*dx$ as $n\to\infty$.
\end{theorem}

\begin{theorem}[Watson--Schwarz identification]\label{thm:intro-dictionary}
If $F_h$ is the inverse of $h$ on the real inverse arc, then
\[
 F_h(x)=x+iP_W(e^{2\pi ix}),\qquad
 \mathcal G(x)=-\frac1{F_h(x)-x}.
\]
\end{theorem}

\begin{proposition}[Horizontal sign]\label{prop:intro-sign}
For $x\notin\frac12\ZZ$,
$\operatorname{sgn}\operatorname{Re}\mathcal G(x)
=\operatorname{sgn}\sin(2\pi x)$.
\end{proposition}

Theorem~\ref{thm:intro-main} is proved in Section~\ref{sec:watson}.
The other two statements are proved in Section~\ref{subsec:real-Schwarz-locus}.

\begin{theorem}[Modular forms and derivatives at a zero]\label{thm:intro-arithmetic}
Let $E_2(\rho)=0$.
\begin{enumerate}[label=(\roman*)]
 \item Evaluation at $\rho$ embeds every finite-level arithmetic modular-function field into $\CC$.
 \item For a quasimodular form with arithmetic meromorphic modular
 coefficients as in Section~\ref{sec:arithmetic-rigidity}, the order at
 $\rho$ equals its $E_2$-adic order.
 \item If $r\in\QQ$ and $E_2(\rho+r)=0$, then $r\in\ZZ$.
 \item One has
 \[
 6j(\rho)(j(\rho)-1728)j''(\rho)
 =(7j(\rho)-6912)j'(\rho)^2,
 \]
 and the algebraic closure in $\CC$ of all finite-level modular-function derivative values equals $\aclC{\Qbar(j(\rho),j'(\rho))}$.
\end{enumerate}
\end{theorem}

Here the arithmetic forms and functions are defined over $\Qbar$.
In (i), every nonzero function is finite and nonzero at $\rho$.
The assertions are proved in Sections~\ref{sec:arithmetic-rigidity},
\ref{sec:ordinary-derivatives}, and \ref{sec:watson}.

For the Legendre parameter $z=\lambda(\tau)$, put
\[
 F_0(z)={}_2F_1\!\left(\frac12,\frac12;1;z\right),
 \qquad
 A(z)=(1-2z)F_0(z)+6z(1-z)F_0'(z).
\]
Then
\[
 E_2(\tau)=F_0(z)A(z),
 \qquad z=\lambda(\tau).
\]
The function $A$ satisfies a second-order Fuchsian equation rationally
gauge-equivalent to the Legendre equation. Its two additional
singularities are apparent, and its Schwarz map is $h$
(Section~\ref{subsec:watson-G}).

Transcendence of the zeros remains open. Eterovi\'c's one-point modular
Schanuel conjecture implies that $\rho,j(\rho),j'(\rho)$ are
algebraically independent. Sections~\ref{sec:conditional-periods}
and \ref{sec:waldschmidt-obstruction} give further conditional
implications and compute the ranks and dimensions in the proposed
application of Waldschmidt's theorem.

\section{Periods, quasi-periods, and \texorpdfstring{$E_2$}{E2}}

Let
\[
 \Lambda=\ZZ\omega_1+\ZZ\omega_2,\quad
 \tau=\frac{\omega_2}{\omega_1}\in\HH.
\]
Let $\zeta_\Lambda$ denote the Weierstrass zeta function. Its quasi-period homomorphism
\[
 H_\Lambda:\Lambda\longrightarrow\CC,
 \quad
 \zeta_\Lambda(z+\omega)=\zeta_\Lambda(z)+H_\Lambda(\omega)
\]
is $\ZZ$-linear. Put
\[
 \eta_1=H_\Lambda(\omega_1),\quad
 \eta_2=H_\Lambda(\omega_2).
\]
With the convention used here, Legendre's relation is
\begin{equation}\label{eq:Legendre}
 \omega_2\eta_1-\omega_1\eta_2=2\pi i.
\end{equation}
Under a homothety $\Lambda\mapsto \alpha\Lambda$, one has
\begin{equation}\label{eq:homogeneity}
 H_{\alpha\Lambda}(\alpha\omega)=\alpha^{-1}H_\Lambda(\omega).
\end{equation}
Thus the ratios $\omega_2/\omega_1$ and $\eta_2/\eta_1$ are unchanged by homothety.

For the normalized lattice $\Lambda_\tau=\ZZ+\ZZ\tau$, one has
\begin{equation}\label{eq:eta-E2-normalized}
 \eta_1(\Lambda_\tau)=G_2(\tau)=\frac{\pi^2}{3}E_2(\tau),
\end{equation}
where
\[
 E_2(\tau)=1-24\sum_{n\geq 1}\sigma_1(n)e^{2\pi i n\tau}.
\]
By homogeneity and Legendre's relation, for the general lattice $\Lambda$,
\begin{equation}\label{eq:eta-E2-general}
 \eta_1=\frac{\pi^2}{3\omega_1}E_2(\tau),
 \quad
 \eta_2=\frac{1}{\omega_1}
 \left(\frac{\pi^2}{3}\tau E_2(\tau)-2\pi i\right).
\end{equation}

Whenever $E_2(\tau)\neq0$, define
\begin{equation}\label{eq:h}
 h(\tau)=\frac{\eta_2}{\eta_1}
 =\tau+\frac{6}{\pi i E_2(\tau)}.
\end{equation}
The function $h$ is equivariant:
\begin{equation}\label{eq:equivariance}
 h(\gamma\tau)=\gamma h(\tau),
 \quad \gamma\in\mathrm{SL}_2(\ZZ),
\end{equation}
where the action on both sides is fractional linear. In projective language, $h(\tau)=\infty$ exactly when $E_2(\tau)=0$. This equivariant formulation is used in \cite{SaberSebbar2012} and the references therein.

\section{Rational quasi-period ratios and zeros of \texorpdfstring{$E_2$}{E2}}

A version of the following equivalence appears in \cite{Kumar2021}. The proof is included for later use.

\begin{proposition}\label{prop:rational-zero}
The following are equivalent.
\begin{enumerate}[label=(\roman*)]
 \item $\eta_1$ and $\eta_2$ are linearly dependent over $\QQ$.
 \item $[\eta_1:\eta_2]\in\PP^1(\QQ)$.
 \item There exists $\gamma\in\mathrm{SL}_2(\ZZ)$ such that $E_2(\gamma\tau)=0$.
\end{enumerate}
\end{proposition}

\begin{proof}
The equivalence of (i) and (ii) is immediate. Suppose
\[
 d\eta_1+c\eta_2=0
\]
with coprime integers $c,d$. Choose $a,b\in\ZZ$ such that $ad-bc=1$ and put
\[
 \gamma=\begin{pmatrix}a&b\\ c&d\end{pmatrix}.
\]
The corresponding change of lattice basis is
\[
 \omega_1'=d\omega_1+c\omega_2,
 \quad
 \omega_2'=b\omega_1+a\omega_2,
\]
so that
\[
 \frac{\omega_2'}{\omega_1'}
 =\frac{a\tau+b}{c\tau+d}=\gamma\tau.
\]
Since the quasi-period map is $\ZZ$-linear,
\[
 \eta_1'=d\eta_1+c\eta_2=0.
\]
Formula \eqref{eq:eta-E2-general} then gives $E_2(\gamma\tau)=0$.

Conversely, if $E_2(\gamma\tau)=0$, then the first quasi-period in the transformed integral basis vanishes. Undoing the integral change of basis yields a nontrivial integral, hence rational, linear relation between $\eta_1$ and $\eta_2$.
\end{proof}

\begin{conjecture}\label{conj:PQ}
Let $\Lambda=\ZZ\omega_1+\ZZ\omega_2$ be a lattice in $\CC$. If
\[
 [\eta_1:\eta_2]\in\PP^1(\QQ),
\]
then
\[
 \frac{\omega_2}{\omega_1}\notin\Qbar.
\]
\end{conjecture}

The hypothesis is equivalent to $\eta_1/\eta_2\in\QQ\cup\{\infty\}$.

\begin{proposition}\label{prop:equiv-zero-conj}
Conjecture \ref{conj:PQ} is equivalent to the assertion that every zero of $E_2$ in $\HH$ is transcendental.
\end{proposition}

\begin{proof}
If $\rho\in\HH$ is a zero of $E_2$, then for $\Lambda_\rho=\ZZ+\ZZ\rho$, formula \eqref{eq:eta-E2-normalized} gives $\eta_1=0$. Thus Conjecture \ref{conj:PQ} implies $\rho\notin\Qbar$.

Conversely, if $\eta_1,\eta_2$ are $\QQ$-linearly dependent, Proposition \ref{prop:rational-zero} gives $\gamma\in\mathrm{SL}_2(\ZZ)$ such that $\rho=\gamma\tau$ is a zero of $E_2$. Fractional linear transformations with rational coefficients preserve algebraicity and transcendence. Hence transcendence of all zeros implies transcendence of $\tau$.
\end{proof}

The zeros of $E_2$ form infinitely many $\mathrm{SL}_2(\ZZ)$-equivalence classes \cite{ElBasraouiSebbar2010}.

\section{Transcendence results at zeros of \texorpdfstring{$E_2$}{E2}}

\begin{proposition}\label{thm:noCM}
No zero of $E_2$ in $\HH$ is imaginary quadratic. Equivalently, no zero of $E_2$ is a CM point.
\end{proposition}

\begin{proof}
Suppose, to the contrary, that $\rho\in\HH$ is imaginary quadratic and
\[
 E_2(\rho)=0.
\]
For the normalized lattice
\[
 \Lambda_\rho=\ZZ+\ZZ\rho,
\]
formula \eqref{eq:eta-E2-normalized} gives
\[
 \eta_1(\Lambda_\rho)=0.
\]
Since $\rho$ is imaginary quadratic, $j(\rho)$ is algebraic. Therefore there exists an elliptic curve $E$ over $\Qbar$ whose complex period lattice is homothetic to $\Lambda_\rho$. Thus, for some $\alpha\in\CC^\times$, a period lattice of a Weierstrass model of $E$ with algebraic invariants $g_2,g_3$ is
\[
 \Lambda'=\alpha\Lambda_\rho.
\]
The nonzero number $\alpha$ is a period of $\Lambda'$. By the homogeneity of the quasi-period map,
\[
 H_{\Lambda'}(\alpha)
 =\alpha^{-1}H_{\Lambda_\rho}(1)
 =\alpha^{-1}\eta_1(\Lambda_\rho)
 =0.
\]
On the other hand, Schneider's theorem \cite{Schneider1937} says that if a Weierstrass elliptic function has algebraic invariants and $\omega\neq0$ is a period, then both $\omega$ and its associated quasi-period $\eta(\omega)$ are transcendental. In particular $\eta(\omega)$ cannot vanish; see \cite[Theorem 18]{Waldschmidt2006}. Applying this to $\omega=\alpha$ gives a contradiction.
\end{proof}

\begin{corollary}\label{cor:notCM-rational}
If $[\eta_1:\eta_2]\in\PP^1(\QQ)$, then $\tau$ is not imaginary quadratic.
\end{corollary}

\begin{proof}
By Proposition \ref{prop:rational-zero}, $\tau$ is modularly equivalent to a zero of $E_2$. The property of being imaginary quadratic is preserved under $\mathrm{SL}_2(\ZZ)$, so Proposition \ref{thm:noCM} applies.
\end{proof}

Nesterenko's theorem \cite{Nesterenko1996} gives the following stronger conclusion.

\begin{proposition}\label{prop:jtrans}
If $E_2(\rho)=0$, then both $j(\rho)$ and $\lambda(\rho)$ are
transcendental. Consequently, for every Weierstrass realization
\[
 y^2=4x^3-g_2x-g_3
\]
of the corresponding homothety class, $g_2$ and $g_3$ cannot both be algebraic.
\end{proposition}

\begin{proof}
Nesterenko's theorem \cite{Nesterenko1996} implies that at least three of
\[
 q=e^{2\pi i\rho},\quad E_2(\rho),\quad E_4(\rho),\quad E_6(\rho)
\]
are algebraically independent over $\Qbar$. Since $E_2(\rho)=0$, it follows that
\[
 q,\ E_4(\rho),\ E_6(\rho)
\]
are algebraically independent. If $j(\rho)$ were algebraic, the identity
\[
 j(\rho)=1728\frac{E_4(\rho)^3}{E_4(\rho)^3-E_6(\rho)^2}
\]
would give the nontrivial algebraic relation
\[
 \bigl(j(\rho)-1728\bigr)E_4(\rho)^3
 -j(\rho)E_6(\rho)^2=0,
\]
a contradiction. Hence $j(\rho)$ is transcendental. The identity
\[
 j(\tau)=256\,
 \frac{(1-\lambda(\tau)+\lambda(\tau)^2)^3}
 {\lambda(\tau)^2(1-\lambda(\tau))^2}
\]
then shows that $\lambda(\rho)$ is transcendental: an algebraic value
of $\lambda(\rho)$ would give an algebraic value of $j(\rho)$. Finally,
since
\[
 j=1728\frac{g_2^3}{g_2^3-27g_3^2},
\]
algebraicity of both $g_2$ and $g_3$ would force $j$ to be algebraic.
\end{proof}

Since CM points have algebraic $j$-invariant,
Proposition~\ref{prop:jtrans} also proves Proposition~\ref{thm:noCM}.
An algebraic zero of $E_2$ would therefore have degree at least $3$.

\section{Modular forms at an \texorpdfstring{$E_2$}{E2}-zero}
\label{sec:arithmetic-rigidity}

Let $\Gamma$ be a congruence subgroup. We call a modular function or
form arithmetic if it is defined over $\Qbar$ in the sense of the
$q$-expansion principle. Choose a number field $F_\Gamma$ over which
the coarse modular curve $X(\Gamma)$ and its finite morphism
$j:X(\Gamma)\to\PP^1$ are defined, and fix the irreducible component
whose complex points contain $\Gamma\backslash\HH$. We write
\[
 K_\Gamma=\Qbar\bigl(X(\Gamma)_{\Qbar}\bigr)
\]
for its arithmetic function field. This also applies when $\Gamma$ is
not neat: stabilizers do not change the coarse function field. Thus
$K_\Gamma/\Qbar(j)$ is finite.

\begin{theorem}[Evaluation on modular curves]
\label{thm:finite-level-genericity}
Let $E_2(\rho)=0$. For every congruence subgroup $\Gamma$, evaluation
at $\rho$ defines an embedding
\[
 K_\Gamma\hookrightarrow\CC.
\]
Equivalently, every nonzero element of $K_\Gamma$ is finite and nonzero
at $\rho$. If $g\in K_\Gamma$ is nonconstant, then $g(\rho)$ is
transcendental and
\[
 \aclC{\Qbar(j(\rho))}
 =\aclC{\Qbar(g(\rho))}.
\]
Moreover, every nonzero arithmetic meromorphic modular form of integral
weight and finite-order character has order zero at $\rho$.
\end{theorem}

\begin{proof}
Let $P_\rho\in X(\Gamma)(\CC)$ be the point represented by $\rho$.
We first show that $P_\rho$ is generic over $\Qbar$.  If it were a
closed point of $X(\Gamma)_{\Qbar}$, then the finite morphism
\[
 j:X(\Gamma)_{\Qbar}\longrightarrow\PP^1_{\Qbar}
\]
would send it to a closed point of $\PP^1_{\Qbar}$.  Because $\Qbar$
is algebraically closed, this would say that $j(\rho)\in\Qbar$, in
contradiction with Proposition~\ref{prop:jtrans}.  Equivalently, the
evaluation homomorphism from the affine coordinate ring of any
nonempty arithmetic chart containing $P_\rho$ has zero kernel.

Now let $0\ne h\in K_\Gamma$.  The zero and pole divisors of $h$ are
defined over $\Qbar$ and are supported on finitely many closed points
of $X(\Gamma)_{\Qbar}$.  Every such point has algebraic $j$-value,
again because $j$ is finite.  Hence $P_\rho$ belongs to neither
divisor.  Thus $h$ is regular and nonzero at $P_\rho$.  Evaluation is
therefore defined on every element of $K_\Gamma$ and commutes with
inversion.  If $h_1\ne h_2$, then $h_1-h_2$ is nonzero and so
$(h_1-h_2)(\rho)\ne0$; consequently the evaluation map is injective.

Suppose next that $g\in K_\Gamma$ is nonconstant.  If
$g(\rho)\in\Qbar$, then $g-g(\rho)$ is a nonzero element of
$K_\Gamma$ vanishing at $P_\rho$, contrary to what was just proved.
Thus $g(\rho)$ is transcendental.  Since $K_\Gamma$ is a function
field of one variable and both $j$ and $g$ are nonconstant, the two
extensions
\[
 K_\Gamma/\Qbar(j),\qquad K_\Gamma/\Qbar(g)
\]
are finite.  Take monic algebraic equations for $g$ over
$\Qbar(j)$ and for $j$ over $\Qbar(g)$.  Their coefficients are
rational functions, all of which can be evaluated at $P_\rho$ by the
first part.  Evaluation therefore shows that $g(\rho)$ is algebraic
over $\Qbar(j(\rho))$ and that $j(\rho)$ is algebraic over
$\Qbar(g(\rho))$.  The two fields consequently have the same
algebraic closure in $\CC$.

Finally, let $f$ be an arithmetic meromorphic modular form of integral
weight $k$ and finite-order character.  Choose $m>0$ such that the
character of $f^m$ is trivial and $km$ is divisible by $12$.  On a
common congruence subgroup the quotient
\[
 u=\frac{f^m}{\Delta^{km/12}}
\]
has weight zero, trivial character, and arithmetic Fourier
expansions, hence belongs to a finite-level arithmetic function field.
The product formula for $\Delta$ shows that $\Delta$ has no zero in
$\HH$, so
\[
 \operatorname{ord}_{\rho}(u)
 =m\operatorname{ord}_{\rho}(f).
\]
The first part gives $\operatorname{ord}_{\rho}(u)=0$, and therefore
$\operatorname{ord}_{\rho}(f)=0$, as required.
\end{proof}

Consider meromorphic quasimodular forms of integral weight $k$ and
depth at most $p$ with an expansion
\begin{equation}\label{eq:KZ-expansion}
 \widetilde f=\sum_{r=0}^{p} f_rE_2^r
\end{equation}
with arithmetic meromorphic modular coefficients $f_r$ of weights
$k-2r$ on one congruence subgroup. The transformation law makes the
coefficients unique.  To see this, choose $N$ such that the subgroup
contains
\[
 \gamma_n=\begin{pmatrix}1&0\\nN&1\end{pmatrix}
 \qquad(n\in\ZZ).
\]
The transformation law of $E_2$ gives
\[
 (E_2|_2\gamma_n)(\tau)
 =E_2(\tau)+X_n(\tau),
 \qquad
 X_n(\tau)=\frac{6nN}{\pi i(nN\tau+1)}.
\]
If $\sum_{r=0}^pf_rE_2^r=0$, applying the weight-$k$ slash operator by
$\gamma_n$ and using the modularity of $f_r$ yields
\[
 \sum_{r=0}^pf_r(\tau)
 \bigl(E_2(\tau)+X_n(\tau)\bigr)^r=0
 \qquad(n\in\ZZ).
\]
For fixed $\tau\in\HH$, the values $X_n(\tau)$ are pairwise distinct.
The polynomial
\[
 \sum_{r=0}^pf_r(\tau)(E_2(\tau)+X)^r
\]
therefore has infinitely many roots and is identically zero.  Its
leading coefficient gives $f_p(\tau)=0$; descending through the
coefficients gives $f_r(\tau)=0$ for every $r$.  Since this holds for
all $\tau$ outside the poles, each $f_r$ vanishes identically. This is
the meromorphic version, at congruence level, of the Kaneko--Zagier
decomposition \cite{KanekoZagier1995}.

\begin{theorem}[Multiplicity at the zeros of $E_2$]
\label{thm:uniform-quasimodular-multiplicity}
For a nonzero form \eqref{eq:KZ-expansion}, put
\[
 \nu_{E_2}(\widetilde f)=\min\{r:f_r\ne0\}.
\]
Then, at every zero $\rho$ of $E_2$,
\[
 \operatorname{ord}_\rho(\widetilde f)
 =\nu_{E_2}(\widetilde f).
\]
In particular, $\widetilde f$ vanishes at one $E_2$-zero if and only if
it vanishes at all of them, equivalently if and only if $E_2$ divides
its expansion; the multiplicity is then independent of the zero.
\end{theorem}

\begin{proof}
Let $s=\nu_{E_2}(\widetilde f)$. Factoring gives
\[
 \widetilde f=E_2^s
 \bigl(f_s+f_{s+1}E_2+\cdots+f_pE_2^{p-s}\bigr).
\]
Write the second factor as $H$.  By
Theorem~\ref{thm:finite-level-genericity}, every coefficient $f_r$ is
regular at $\rho$, and $f_s(\rho)\ne0$.  Since $E_2(\rho)=0$, it
follows that
\[
 H(\rho)=f_s(\rho)\ne0.
\]
Thus $H$ is a unit in the local ring of holomorphic germs at $\rho$,
and
\[
 \operatorname{ord}_{\rho}(\widetilde f)
 =s\operatorname{ord}_{\rho}(E_2).
\]

It remains only to determine the order of $E_2$.  The Ramanujan identity
\[
 D E_2=\frac{E_2^2-E_4}{12},
 \qquad D=\frac1{2\pi i}\frac d{d\tau},
\]
gives
\[
 D E_2(\rho)=-\frac{E_4(\rho)}{12}.
\]
Proposition~\ref{prop:jtrans} gives $E_4(\rho)\ne0$. Hence
$E_2'(\rho)=2\pi iD E_2(\rho)\ne0$, so $E_2$ has a simple zero at
$\rho$.  Therefore $\operatorname{ord}_{\rho}(\widetilde f)=s$.
The final equivalences follow because $s\ge1$ holds exactly when the
coefficient $f_0$ vanishes, that is, exactly when the expansion is
divisible by $E_2$.
\end{proof}

\section{Ordinary derivatives and modular Schanuel}
\label{sec:ordinary-derivatives}

Primes denote ordinary $\tau$-derivatives. We use the Ramanujan
equations and the identities
\[
 D\log\Delta=E_2,\qquad
 \frac{Dj}{j}=-\frac{E_6}{E_4}.
\]

\begin{proposition}
\label{prop:j-second-derivative}
On the open set where $j(j-1728)Dj\ne0$,
\begin{equation}\label{eq:E2-from-j-derivatives}
 E_2=6\frac{D^2j}{Dj}-4\frac{Dj}{j}
                  -3\frac{Dj}{j-1728}.
\end{equation}
Consequently every zero $\rho$ of $E_2$ satisfies
\begin{equation}\label{eq:j-second-derivative}
 6j(\rho)(j(\rho)-1728)j''(\rho)
 =(7j(\rho)-6912)j'(\rho)^2.
\end{equation}
\end{proposition}

\begin{proof}
Put $r=E_6/E_4$.  From
\[
 j=\frac{E_4^3}{\Delta},\qquad D\log\Delta=E_2,
 \qquad DE_4=\frac{E_2E_4-E_6}{3},
\]
we obtain
\[
 \frac{Dj}{j}=3\frac{DE_4}{E_4}-\frac{D\Delta}{\Delta}
 =-\frac{E_6}{E_4}=-r.
\]
The identities $j=E_4^3/\Delta$ and
$j-1728=E_6^2/\Delta$ also give
\[
 E_4=\frac{jr^2}{j-1728},\qquad E_6=rE_4,
\]
because $(j-1728)/j=r^2/E_4$.  Consequently
\[
 \Delta=\frac{E_4^3}{j}
 =\frac{j^2r^6}{(j-1728)^3}
 =\frac{(Dj)^6}{j^4(j-1728)^3}.
\]
Taking the normalized logarithmic derivative of this last identity
yields
\[
 E_2=D\log\Delta
 =6\frac{D^2j}{Dj}-4\frac{Dj}{j}
  -3\frac{Dj}{j-1728},
\]
which proves \eqref{eq:E2-from-j-derivatives} on the stated open set.

Let now $E_2(\rho)=0$.  The algebraic independence used in
Proposition~\ref{prop:jtrans} gives
$E_4(\rho)E_6(\rho)\ne0$.  Moreover $j(\rho)$ is transcendental, so
$j(\rho)\notin\{0,1728\}$, and
\[
 Dj(\rho)=-j(\rho)\frac{E_6(\rho)}{E_4(\rho)}\ne0.
\]
Thus the preceding identity may be evaluated at $\rho$.  Clearing its
denominators gives
\[
 6j(\rho)(j(\rho)-1728)D^2j(\rho)
 =(7j(\rho)-6912)Dj(\rho)^2.
\]
Finally $Dj=j'/(2\pi i)$ and
$D^2j=j''/(2\pi i)^2$; the common factor $(2\pi i)^{-2}$ cancels and
gives \eqref{eq:j-second-derivative}.
\end{proof}

Higher derivatives are obtained from the Ramanujan differential
equations. The identity \eqref{eq:j-second-derivative}, which holds
only at the zeros of $E_2$, cannot itself be differentiated.

\begin{proposition}[All higher $j$-derivatives]
\label{prop:j-higher-derivatives}
For every $n\ge1$ there is a rational function $R_n\in\QQ(J)$,
independent of $\rho$, such that
\begin{equation}\label{eq:j-higher-derivatives}
 j^{(n)}(\rho)=j'(\rho)^nR_n(j(\rho))
\end{equation}
whenever $E_2(\rho)=0$. In particular, every $j^{(n)}(\rho)$ is
algebraic over $\QQ(j(\rho),j'(\rho))$, and
\[
 R_1(J)=1,\qquad
 R_2(J)=\frac{7J-6912}{6J(J-1728)}.
\]
\end{proposition}

\begin{proof}
In the field of meromorphic functions put
\[
 J=j,\qquad r=\frac{E_6}{E_4},\qquad U=\frac{E_2}{r},
\]
and set
\[
 a(J)=\frac{J+3456}{6(J-1728)},\qquad
 b(J)=\frac{J}{J-1728}.
\]
The Ramanujan equations
\[
 DE_2=\frac{E_2^2-E_4}{12},\qquad
 DE_4=\frac{E_2E_4-E_6}{3},\qquad
 DE_6=\frac{E_2E_6-E_4^2}{2}
\]
give the differential system:
\begin{equation}\label{eq:JU-system-Annalen}
\begin{aligned}
 DJ&=-Jr,\\
 Dr&=r^2\bigl(U/6-a(J)\bigr),\\
 DU&=r\bigl(-U^2/12+a(J)U-b(J)/12\bigr).
\end{aligned}
\end{equation}
Indeed, the first equation is $Dj/j=-E_6/E_4$.  For the second one,
differentiate $r=E_6/E_4$ and substitute the Ramanujan equations:
\[
\begin{aligned}
 Dr
 &=\frac{(DE_6)E_4-E_6(DE_4)}{E_4^2}\\
 &=\frac{E_2r}{6}-\frac{E_4}{2}+\frac{r^2}{3}.
\end{aligned}
\]
Now use
\[
 E_2=rU,\qquad E_4=\frac{Jr^2}{J-1728}.
\]
This gives $Dr=r^2(U/6-a(J))$.  Similarly,
\[
\begin{aligned}
 DU
 &=\frac{DE_2}{r}-\frac{E_2Dr}{r^2}\\
 &=r\left(\frac{U^2}{12}-\frac{b(J)}{12}
           -U\bigl(U/6-a(J)\bigr)\right)\\
 &=r\bigl(-U^2/12+a(J)U-b(J)/12\bigr),
\end{aligned}
\]
which proves the system.

Starting with $A_0(J,U)=J$, define
\begin{equation}\label{eq:A-recurrence-Annalen}
\begin{aligned}
 A_{n+1}={}&-J\partial_JA_n
 +\bigl(-U^2/12+a(J)U-b(J)/12\bigr)\partial_UA_n\\
 &+n\bigl(U/6-a(J)\bigr)A_n.
\end{aligned}
\end{equation}
We claim that
\begin{equation}\label{eq:Dn-j-An}
 D^nj=r^nA_n(J,U)\qquad(n\ge0).
\end{equation}
The identity holds for $n=0$.  If it holds for $n$, the chain rule and
\eqref{eq:JU-system-Annalen} give
\[
\begin{aligned}
 D^{n+1}j
 &=r^n\bigl((DJ)\partial_JA_n+(DU)\partial_UA_n\bigr)
   +nr^{n-1}(Dr)A_n\\
 &=r^{n+1}A_{n+1}(J,U),
\end{aligned}
\]
which is exactly recurrence \eqref{eq:A-recurrence-Annalen}.  This
proves \eqref{eq:Dn-j-An} by induction.

At an $E_2$-zero, $U=0$ because $r(\rho)\ne0$, and
$2\pi ir=-j'/J$.  Multiplying \eqref{eq:Dn-j-An} by
$(2\pi i)^n$ therefore gives
\[
 j^{(n)}(\rho)
 =j'(\rho)^n\left(-\frac1{J}\right)^nA_n(J,0)
 \quad\text{at }J=j(\rho).
\]
Thus $R_n(J)=(-1/J)^nA_n(J,0)$ proves
\eqref{eq:j-higher-derivatives}.  To identify the first two
terms, the recurrence gives
\[
 A_1=-J,\qquad
 A_2=J\bigl(1-U/6+a(J)\bigr).
\]
Hence $R_1=1$ and
\[
 R_2(J)=\frac{1+a(J)}{J}
 =\frac{7J-6912}{6J(J-1728)}.
\]
In particular every $j^{(n)}(\rho)$ belongs to
$\QQ(j(\rho),j'(\rho))$, and the asserted algebraicity follows.
\end{proof}

\medskip
The finite-level fields form a directed system: if $\Gamma_1$ and
$\Gamma_2$ are congruence subgroups, both arithmetic function fields
embed in the field for $\Gamma_1\cap\Gamma_2$. Thus any finite
collection of modular functions and their derivatives may be treated
after passage to one common congruence subgroup.

\begin{theorem}
\label{thm:finite-level-derivative-tower}
Let $\mathcal F^\partial_{\rm fin}(\rho)$ be the field generated over
$\Qbar$ by all values $g^{(n)}(\rho)$, where $n\ge0$ and $g$ ranges
over the arithmetic meromorphic modular functions on congruence
subgroups. Then
\[
 \aclC{\mathcal F^\partial_{\rm fin}(\rho)}
 =\aclC{\Qbar(j(\rho),j'(\rho))}.
\]
\end{theorem}

\begin{proof}
It is enough to treat finitely many generators at a time.  After
intersecting their congruence subgroups, they all belong to one
arithmetic function field $K_\Gamma$.  The extension
$K_\Gamma/\Qbar(j)$ is finite and separable, since the ground field has
characteristic zero.  Hence the derivation $d/dj$ of $\Qbar(j)$ has a
unique extension
\[
 \delta:K_\Gamma\longrightarrow K_\Gamma.
\]
It is $\Qbar$-linear and, on analytic functions away from ramification
points, satisfies
\begin{equation}\label{eq:delta-chain-rule}
 \frac{dg}{d\tau}=j'\delta g.
\end{equation}
Both sides are meromorphic, so this identity holds identically and can
be used without excluding the ramification points.

Iterating \eqref{eq:delta-chain-rule} expresses $g^{(n)}$ as a
universal polynomial with integer coefficients in
\[
 j',j'',\ldots,j^{(n)}\quad\text{and}\quad
 \delta g,\delta^2g,\ldots,\delta^ng.
\]
For example,
\[
 g'=j'\delta g,\qquad
 g''=j''\delta g+(j')^2\delta^2g,
\]
and the general assertion follows by differentiating once more and
using \eqref{eq:delta-chain-rule}.  Each $\delta^kg$ belongs to
$K_\Gamma$.  Since $K_\Gamma$ is finite over $\Qbar(j)$,
Theorem~\ref{thm:finite-level-genericity} shows that
$(\delta^kg)(\rho)$ is algebraic over $\Qbar(j(\rho))$; the same
theorem also ensures that the evaluation is defined even when the
chosen representative has poles elsewhere on the modular curve.

Proposition~\ref{prop:j-higher-derivatives} places every
$j^{(r)}(\rho)$ in $\QQ(j(\rho),j'(\rho))$.  The displayed polynomial
expression therefore proves that every generator $g^{(n)}(\rho)$ of
$\mathcal F^\partial_{\rm fin}(\rho)$ is algebraic over
$\Qbar(j(\rho),j'(\rho))$.  This gives one inclusion between the
algebraic closures.  For the reverse inclusion, take $g=j$ and
$n=0,1$: both $j(\rho)$ and $j'(\rho)$ occur among the permitted
generators.  The two algebraic closures are consequently equal.
\end{proof}

\medskip
\begin{conjecture}[Eterovi\'c's one-point modular Schanuel conjecture]
\label{conj:Eterovic-one-point}
For every non-CM $\tau\in\HH$,
\begin{equation}\label{eq:Eterovic-one-point}
 \trdeg_\QQ\QQ\bigl(\tau,j(\tau),j'(\tau),j''(\tau)\bigr)\ge3.
\end{equation}
\end{conjecture}

This is the one-point case of
\cite[Conjecture~2.3]{Eterovic2025}. It is a conjecture about values;
the functional Ax--Schanuel theorem \cite{PilaTsimerman2016} does not
give this inequality at a prescribed point.

\begin{theorem}[Conditional algebraic independence]
\label{thm:Eterovic-consequence}
Assume Conjecture~\ref{conj:Eterovic-one-point}. If $E_2(\rho)=0$, then
\[
 \rho,\qquad j(\rho),\qquad j'(\rho)
\]
are algebraically independent over $\Qbar$. In particular, $\rho$ is
transcendental.
\end{theorem}

\begin{proof}
Proposition~\ref{prop:jtrans} shows that $\rho$ is not CM, because a
CM point has algebraic $j$-invariant whereas $j(\rho)$ is
transcendental.  Thus Conjecture~\ref{conj:Eterovic-one-point} applies
at $\rho$.  Moreover $j(\rho)\notin\{0,1728\}$, and
equation~\eqref{eq:j-second-derivative} may be solved for $j''(\rho)$:
\[
 j''(\rho)=
 \frac{7j(\rho)-6912}{6j(\rho)(j(\rho)-1728)}j'(\rho)^2.
\]
Consequently
\[
 \QQ(\rho,j(\rho),j'(\rho),j''(\rho))
 =\QQ(\rho,j(\rho),j'(\rho)).
\]
The conjecture makes the field on the left have transcendence degree
at least three over $\QQ$.  The field on the right is generated by
three numbers, so its transcendence degree is at most three.  Equality
therefore holds, which is equivalent to the algebraic independence of
$\rho,j(\rho),j'(\rho)$ over $\QQ$.  Algebraic independence is
unchanged after the algebraic extension $\QQ\subset\Qbar$, giving the
stated conclusion.
\end{proof}

\section{Watson's expansion}
\label{sec:watson}

Watson's inversion formula \cite[\S\S4--6, pp.~137--149]{Watson1952}
parametrizes the zeros of $E_2$, modulo translation, by rational cusps.
Let $m/n\in\QQ$ be reduced, with $n>0$, and choose the centered
integer $m'$ satisfying
\begin{equation}\label{eq:Watson-inverse}
 mm'\equiv -1\pmod n,
 \quad
 -\frac n2<m'\leq\frac n2,
\end{equation}
with $m'=0$ for $m/n=0/1$. Watson constructs real coefficients $c_r$ for which the series below is absolutely convergent and writes the transformed point in the fundamental region as
\begin{equation}\label{eq:Watson-T}
 T_{m/n}
 =\frac{m'}n+\frac{6i}{\pi}
 +i\sum_{r\geq1}(-1)^{r-1}c_r
 e^{2\pi i r m'/n}.
\end{equation}
His Lagrange expansion gives
\begin{equation}\label{eq:Watson-mobius}
 \left(T_{m/n}-\frac{m'}n\right)
 \left(\rho_{m/n}-\frac mn\right)
 =-\frac1{n^2},
\end{equation}
where $\rho_{m/n}$ denotes the zero associated with the cusp $m/n$. Watson also proves majorant estimates implying
\begin{equation}\label{eq:Watson-cr-bound}
 \sum_{r\ge1}|c_r|<2.83\times10^{-4}.
\end{equation}

Define
\begin{equation}\label{eq:Watson-P}
 P_W(z)=\frac6\pi+
 \sum_{r\ge1}(-1)^{r-1}c_r z^r,
 \quad |z|\leq1,
\end{equation}
and the one-periodic function
\begin{equation}\label{eq:Watson-profile-function}
 \mathcal G(x)=\frac{i}{P_W(e^{2\pi i x})},
 \quad x\in\RR.
\end{equation}
The bound \eqref{eq:Watson-cr-bound} shows that $P_W$ has no zero on the closed unit disk, since $6/\pi>1.9$. Absolute convergence gives continuity on the closed disk, so $\mathcal G$ is continuous and one-periodic.

By \cite[\S2]{ImamogluJermannToth2014} and
\cite[Theorems~5.1 and 6.2]{WoodYoung2014}, $h$ maps the real inverse
arc in the standard fundamental domain bijectively onto $[-1/2,1/2]$
and has nonzero derivative there. The holomorphic inverse-function
theorem gives local real-analytic inverses, which agree by uniqueness.
Write $F_h$ for this inverse. The
translation law $h(\tau+1)=h(\tau)+1$ extends them to all of $\RR$ by
$F_h(x+1)=F_h(x)+1$, including real-analytic continuation across the
endpoints.

For rational $x=m'/n$ in the centered interval, Watson's transformed
point in \eqref{eq:Watson-T} satisfies $h(T_{m/n})=x$. Uniqueness on
the inverse arc gives
\[
 F_h(x)=x+iP_W(e^{2\pi ix}).
\]
Every rational residue occurs as an inverse residue in
\eqref{eq:Watson-inverse}. Both sides are continuous, so density and
translation extend this equality to every real $x$. It follows that
\[
 \mathcal G(x)=-\frac{1}{F_h(x)-x}.
\]
Since $F_h(x)\in\HH$, the denominator is nonzero. Hence
$\mathcal G$ is real analytic.

Watson's formula now takes the following form.

\begin{proposition}\label{prop:Watson-profile}
For every reduced fraction $m/n$, with $m'$ as in \eqref{eq:Watson-inverse},
\begin{equation}\label{eq:Watson-profile}
 \rho_{m/n}
 =\frac mn+\frac1{n^2}
 \mathcal G\!\left(\frac{m'}n\right).
\end{equation}
Consequently
\begin{equation}\label{eq:Watson-field-equivalence}
 \rho_{m/n}\in\Qbar
 \quad\Longleftrightarrow\quad
 \mathcal G(m'/n)\in\Qbar,
\end{equation}
and when these equivalent conditions hold,
\begin{equation}\label{eq:Watson-degree-equality}
 [\QQ(\rho_{m/n}):\QQ]
 =
 [\QQ(\mathcal G(m'/n)):\QQ].
\end{equation}
Moreover, as $m$ runs through the reduced residue classes modulo $n$, $m'$ does also. Hence Conjecture~\ref{conj:PQ} is equivalent to
\begin{equation}\label{eq:G-rational-trans}
 \mathcal G(x)\notin\Qbar
 \quad\text{for every }x\in\QQ.
\end{equation}
\end{proposition}

\begin{proof}
Equation \eqref{eq:Watson-T} says
\[
 T_{m/n}-\frac{m'}n
 =iP_W(e^{2\pi i m'/n}).
\]
Substitution in \eqref{eq:Watson-mobius} gives
\[
 \rho_{m/n}-\frac mn
 =-\frac{1}{n^2iP_W(e^{2\pi i m'/n})}
 =\frac1{n^2}\mathcal G\!\left(\frac{m'}n\right),
\]
which is \eqref{eq:Watson-profile}.  Since $m/n$ and $n$ are
rational, this identity and its rearrangement
\[
 \mathcal G\!\left(\frac{m'}n\right)
 =n^2\left(\rho_{m/n}-\frac mn\right)
\]
show that
\[
 \QQ(\rho_{m/n})
 =\QQ\!\left(\mathcal G(m'/n)\right).
\]
This proves both the algebraicity equivalence and, when the values are
algebraic, the equality of degrees.

Finally $m'\equiv-m^{-1}\pmod n$.  Inversion followed by
multiplication by $-1$ is a permutation of
$(\ZZ/n\ZZ)^\times$, so $m'$ runs through all reduced residue classes
when $m$ does.  Conversely every rational point of $\RR/\ZZ$ is
represented by such a reduced residue class, with $0$ supplied by the
cusp $0/1$.  The periodicity of $\mathcal G$ now proves the final
equivalence with Conjecture~\ref{conj:PQ}.
\end{proof}

The sharp estimate of Imamo\u{g}lu, Jermann and T\'oth \cite[Theorem~1.1, formula~(5)]{ImamogluJermannToth2014} becomes, in this notation,
\begin{equation}\label{eq:G-IJT-bound}
 \left|\mathcal G(x)-\frac{\pi i}{6}\right|
 \leq 7.77\times10^{-5}
 \quad(x\in\RR).
\end{equation}
The estimate holds first for rational $x$ and extends to all real
$x$ by continuity.

\begin{proposition}\label{prop:Watson-mean}
One has
\begin{equation}\label{eq:Watson-mean}
 \int_0^1\mathcal G(x)\,dx=\frac{\pi i}{6}.
\end{equation}
In particular, if $\nu(x)=\operatorname{Im}\mathcal G(x)$, then
\begin{equation}\label{eq:h-mean}
 \int_0^1 \nu(x)\,dx=\frac\pi6,
 \quad
 \int_0^1\operatorname{Re}\mathcal G(x)\,dx=0.
\end{equation}
\end{proposition}

\begin{proof}
The series defining $P_W$ converges uniformly on the closed unit disk
by \eqref{eq:Watson-cr-bound}.  Moreover,
\[
 |P_W(z)|\ge \frac6\pi-\sum_{r\ge1}|c_r|>0
 \qquad(|z|\le1).
\]
Hence $F(z)=i/P_W(z)$ is holomorphic for $|z|<1$ and continuous for
$|z|\le1$.  For $0<r<1$, Cauchy's mean-value formula gives
\[
 \int_0^1F(re^{2\pi ix})\,dx=F(0)
 =\frac{i}{P_W(0)}=\frac{\pi i}{6}.
\]
Because $F$ is uniformly continuous on the closed disk,
$F(re^{2\pi ix})$ converges uniformly to $F(e^{2\pi ix})$ as
$r\to1^-$.  We may therefore pass to the limit under the integral.
Since $F(e^{2\pi ix})=\mathcal G(x)$, this proves
\eqref{eq:Watson-mean}.  Taking real and imaginary parts gives
\eqref{eq:h-mean}.
\end{proof}

The height count follows from the cusp parametrization and the
following lattice-point estimate.

\begin{lemma}\label{lem:primitive-lattice-count}
Let $\mathcal R\subset\mathbb R^2$ be bounded with piecewise $C^1$ boundary and suppose that $0\notin\mathcal R$. Then, as $T\to\infty$,
\[
 \#\bigl(T\mathcal R\cap\mathbb Z^2_{\rm prim}\bigr)
 =\frac{\operatorname{area}(\mathcal R)}{\zeta(2)}T^2
 +O_{\mathcal R}(T\log T).
\]
\end{lemma}

\begin{proof}
The ordinary lattice-point estimate for a bounded planar region with piecewise $C^1$ boundary is
\[
 \#(t\mathcal R\cap\mathbb Z^2)
 =\operatorname{area}(\mathcal R)t^2+O_{\mathcal R}(t+1).
\]
Apply M\"obius inversion to the gcd of the two coordinates. Since every nonzero lattice point of $T\mathcal R$ has gcd $O_{\mathcal R}(T)$,
\[
 \#\bigl(T\mathcal R\cap\mathbb Z^2_{\rm prim}\bigr)
 =\sum_{d\ll_{\mathcal R}T}\mu(d)
 \#\bigl((T/d)\mathcal R\cap\mathbb Z^2\bigr).
\]
The main terms give
\[
 \operatorname{area}(\mathcal R)T^2
 \sum_{d\ll T}\frac{\mu(d)}{d^2}
 =\frac{\operatorname{area}(\mathcal R)}{\zeta(2)}T^2+O_{\mathcal R}(T),
\]
while the boundary errors contribute
\[
 O_{\mathcal R}\!\left(\sum_{d\ll T}(T/d+1)\right)
 =O_{\mathcal R}(T\log T).
\]
\end{proof}

\begin{theorem}\label{thm:Watson-counting}
Let
\[
 N(Y)=\#\left\{\rho:\ E_2(\rho)=0,\ 0\leq\operatorname{Re}\rho<1,
 \ \operatorname{Im}\rho\geq Y\right\}.
\]
Then, as $Y\to0^+$,
\begin{equation}\label{eq:Watson-counting}
 N(Y)=\frac{1}{2\pi Y}
 +O\!\left(Y^{-1/2}\log\frac1Y\right).
\end{equation}
In particular,
\[
 N(Y)\sim\frac1{2\pi Y}.
\]
\end{theorem}

\begin{proof}
By the one-zero-per-cusp theorem of Imamo\u{g}lu, Jermann and T\'oth
\cite[Theorem~1.1]{ImamogluJermannToth2014}, the indexing set is exactly
\[
 \{0/1\}\ \sqcup\
 \{m/n:n\ge2,\ 1\le m<n,\ (m,n)=1\}.
\]
From \eqref{eq:Watson-profile},
\[
 \operatorname{Im}\rho_{m/n}
 =\frac1{n^2}\nu\!\left(\frac{m'}n\right).
\]
For the counting argument, let $\widetilde m'$ denote the representative of $m'$ modulo $n$ with
\[
 0<\widetilde m'<n.
\]
Since $\nu$ is one-periodic,
\[
 \nu\!\left(\frac{m'}n\right)
 =\nu\!\left(\frac{\widetilde m'}n\right),
\]
and the map $m\mapsto\widetilde m'$ permutes the reduced residue classes modulo $n$. Hence, exactly,
\begin{equation}\label{eq:N-endpoint-decomposition}
 N(Y)=\mathbf 1_{\{\operatorname{Im}\rho_{0/1}\geq Y\}}
 +\#\left\{(\widetilde m',n)\in\ZZ^2_{\rm prim}:n\ge2,
 \ 0<\widetilde m'<n,
 \ n^2Y\leq\nu(\widetilde m'/n)\right\}.
\end{equation}
Put $T=Y^{-1/2}$ and consider the bounded region
\[
 \mathcal R
 =\left\{(u,v)\in\RR^2:
 0<u<v,\quad
 v\leq\sqrt{\nu(u/v)}\right\}.
\]
The uniform bound \eqref{eq:G-IJT-bound} implies $\nu(x)>0$. More
explicitly, the closure of $\mathcal R$ has boundary formed by the two
line segments
\[
 \{(0,v):0\leq v\leq\sqrt{\nu(0)}\},\qquad
 \{(v,v):0\leq v\leq\sqrt{\nu(1)}\},
\]
and the real-analytic arc
\[
 \gamma(x)=\bigl(x\sqrt{\nu(x)},\sqrt{\nu(x)}\bigr),
 \qquad 0\leq x\leq1.
\]
 Thus the boundary is
piecewise $C^1$ and has finite length, including at the corner at the
origin. The condition $\operatorname{Im}\rho_{m/n}\geq Y$ is precisely
that $(\widetilde m',n)$ lie in $T\mathcal R$. Equation
\eqref{eq:N-endpoint-decomposition} differs from this primitive lattice
count only by the single endpoint contribution, which is $O(1)$ (and
equals one for all sufficiently small $Y$).

With $u=xv$, the Jacobian is $v$, and Proposition~\ref{prop:Watson-mean} gives
\[
 \operatorname{area}(\mathcal R)
 =\int_0^1\int_0^{\sqrt{\nu(x)}}v\,dv\,dx
 =\frac12\int_0^1\nu(x)\,dx
 =\frac\pi{12}.
\]
Lemma~\ref{lem:primitive-lattice-count} gives
\[
 \#\bigl(T\mathcal R\cap\ZZ^2_{\rm prim}\bigr)
 =\frac{1}{\zeta(2)}\operatorname{area}(\mathcal R)T^2
 +O(T\log T).
\]
Since $1/\zeta(2)=6/\pi^2$, the main coefficient is
\[
 \frac6{\pi^2}\cdot\frac\pi{12}=\frac1{2\pi},
\]
which proves the result.
\end{proof}

\begin{corollary}\label{cor:Watson-equidistribution}
For $n\ge2$ define the probability measure
\[
 \mu_n=\frac1{\varphi(n)}
 \sum_{\substack{0\leq m<n\\(m,n)=1}}
 \delta_{\,n^2(\rho_{m/n}-m/n)}.
\]
Then, for every bounded continuous function $f:\CC\to\CC$,
\begin{equation}\label{eq:Watson-equidistribution}
 \lim_{n\to\infty}\int_{\CC} f(z)\,d\mu_n(z)
 =\int_0^1 f\bigl(\mathcal G(x)\bigr)\,dx.
\end{equation}
Equivalently, $\mu_n$ converges weakly to the pushforward of normalized Lebesgue measure on $[0,1]$ under $\mathcal G$. In particular,
\begin{equation}\label{eq:Watson-average-height}
 \frac1{\varphi(n)}
 \sum_{\substack{0\leq m<n\\(m,n)=1}}
 n^2\operatorname{Im}\rho_{m/n}
 \longrightarrow\frac\pi6.
\end{equation}
\end{corollary}

\begin{proof}
By Proposition~\ref{prop:Watson-profile},
\[
 n^2\left(\rho_{m/n}-\frac mn\right)
 =\mathcal G\!\left(\frac{m'}n\right).
\]
For fixed $n$, the map $m\mapsto m'$ defined by $mm'\equiv-1\pmod n$ permutes the reduced residue classes modulo $n$. Hence, for every bounded continuous $f$,
\[
 \int_{\CC}f(z)\,d\mu_n(z)
 =\frac1{\varphi(n)}
 \sum_{\substack{0\leq a<n\\(a,n)=1}}
 f\!\left(\mathcal G\!\left(\frac an\right)\right).
\]
The reduced residue classes are equidistributed modulo one as $n\to\infty$. Indeed, for every fixed nonzero integer $k$, the Ramanujan sum
\[
 c_n(k)=\sum_{\substack{0\leq a<n\\(a,n)=1}}e^{2\pi i k a/n}
\]
satisfies the quantitative bound
\[
 \left|\frac{c_n(k)}{\varphi(n)}\right|
 \leq\frac{\sigma_1((n,k))}{\varphi(n)}
 \leq\frac{\sigma_1(|k|)}{\varphi(n)}\longrightarrow0.
\]
Weyl's criterion applies. Since $f\circ\mathcal G$ is continuous and one-periodic, the preceding average therefore tends to
\[
 \int_0^1 f\bigl(\mathcal G(x)\bigr)\,dx,
\]
which proves \eqref{eq:Watson-equidistribution}.

For the final assertion, apply the same argument to $x\mapsto\operatorname{Im}\mathcal G(x)$. Using Proposition~\ref{prop:Watson-mean},
\[
 \lim_{n\to\infty}
 \frac1{\varphi(n)}
 \sum_{\substack{0\leq m<n\\(m,n)=1}}
 n^2\operatorname{Im}\rho_{m/n}
 =\int_0^1\operatorname{Im}\mathcal G(x)\,dx
 =\frac\pi6.
\]
\end{proof}

\subsection{Rational separation and rational real parts}

\begin{lemma}\label{lem:rational-translate-modular-form}
Let $r=A/N\in\mathbb Q\setminus\mathbb Z$ be in lowest terms, with $N\ge2$, and put
\[
 F_r(\tau)=E_2(\tau+r)-E_2(\tau).
\]
Then $F_r$ is a nonzero arithmetic holomorphic modular form of weight
$2$ on $\Gamma(N^2)$, defined over the cyclotomic field
$\mathbb Q(\zeta_{N^2})$.
\end{lemma}

\begin{proof}
Let
\[
 \alpha_r=\begin{pmatrix}1&r\\0&1\end{pmatrix}
\]
and let $\gamma=\left(\begin{smallmatrix}a&b\\c&d\end{smallmatrix}\right)\in\Gamma(N^2)$. Then
\[
 \delta=\alpha_r\gamma\alpha_r^{-1}
 =\begin{pmatrix}
 a+rc & b+r(d-a)-r^2c\\
 c & d-rc
 \end{pmatrix}
 \in\mathrm{SL}_2(\mathbb Z),
\]
because $b,c,d-a$ are divisible by $N^2$. Since $\alpha_r\gamma=\delta\alpha_r$ and
\[
 c(\tau+r)+d-rc=c\tau+d,
\]
the anomalous terms in the transformation law of $E_2$ cancel, giving
\begin{equation}\label{eq:rational-translate-modular}
 F_r(\gamma\tau)=(c\tau+d)^2F_r(\tau).
\end{equation}

For holomorphy at the cusps, put
\[
 U_r(\tau)=\frac{\Delta(\tau+r)}{\Delta(\tau)}.
\]
Since $(2\pi i)^{-1}(\log\Delta)'=E_2$,
\begin{equation}\label{eq:Fr-logarithmic-derivative}
 F_r(\tau)=\frac{1}{2\pi i}\frac{d}{d\tau}\log U_r(\tau).
\end{equation}
The same conjugation calculation shows that $U_r$ is invariant under $\Gamma(N^2)$. Moreover, both $\Delta(\tau)$ and $\Delta(\tau+r)=\Delta|_{12}\alpha_r$ restrict to cusp forms on $\Gamma(N^2)$, so $U_r$ extends to a modular unit on the compact modular curve $X(N^2)$, with divisor supported at the cusps. At a cusp $\mathfrak a$, with local parameter $q_{\mathfrak a}$ and width $w_{\mathfrak a}$, write
\[
 U_r=q_{\mathfrak a}^{m}u(q_{\mathfrak a}),\qquad u(0)\ne0.
\]
Then \eqref{eq:Fr-logarithmic-derivative} gives
\[
 (F_r|_2\sigma)(\tau)
 =\frac1{w_{\mathfrak a}}
 \left(m+q_{\mathfrak a}\frac{u'(q_{\mathfrak a})}{u(q_{\mathfrak a})}\right),
\]
Indeed, $q_{\mathfrak a}=e^{2\pi i\tau/w_{\mathfrak a}}$, so the
chain rule gives $(2\pi i)^{-1}d/d\tau=w_{\mathfrak a}^{-1}
q_{\mathfrak a}d/dq_{\mathfrak a}$ after applying the weight-two
slash operator. This also shows directly that no negative power of the
local parameter occurs.  Thus $F_r$ is holomorphic at every cusp. Its
Fourier expansion at infinity is
\begin{equation}\label{eq:rational-translate-qexp}
 F_r(\tau)
 =-24\sum_{n\ge1}\sigma_1(n)
 \bigl(e^{2\pi i nr}-1\bigr)q^n,
 \qquad q=e^{2\pi i\tau},
\end{equation}
so its coefficients lie in
$\mathbb Q(\zeta_N)\subset\mathbb Q(\zeta_{N^2})$, and the coefficient
of $q$ shows that it is nonzero.
\end{proof}

\begin{theorem}\label{thm:rational-real-parts}
Let $\rho\in\mathbb H$ satisfy $E_2(\rho)=0$. If $\operatorname{Re}\rho\in\mathbb Q$, then
\[
 2\operatorname{Re}\rho\in\mathbb Z.
\]
Equivalently, modulo integer translation, every zero of $E_2$ with rational real part lies on one of the two lines
\[
 \operatorname{Re}\tau=0,
 \qquad
 \operatorname{Re}\tau=\frac12.
\]
Consequently, in the strip $0\le\operatorname{Re}\tau<1$, the only zeros with rational real part are $\rho_0$ and $\rho_{1/2}$ defined below.
\end{theorem}

\begin{proof}
Write $x=\operatorname{Re}\rho\in\mathbb Q$. Since $E_2$ has real Fourier coefficients,
\[
 E_2(-\overline\rho)=\overline{E_2(\rho)}=0.
\]
But $-\overline\rho=\rho-2x$. Put $r=-2x$. If $r\notin\mathbb Z$, Lemma~\ref{lem:rational-translate-modular-form} applied to this rational translation gives
a nonzero arithmetic modular form $F_r$ with $F_r(\rho)=0$. This
contradicts Theorem~\ref{thm:finite-level-genericity}. Hence
$r\in\mathbb Z$, which is equivalent to $2x\in\mathbb Z$.

The final assertion follows from the uniqueness of the zero on each reflection-symmetric line, proved in \cite{ElBasraouiSebbar2010} and also encoded in Watson's cusp parametrization \cite{Watson1952}.
\end{proof}

\begin{corollary}\label{cor:no-rational-separated-zeros}
If $E_2(\rho)=E_2(\rho+r)=0$ with $r\in\mathbb Q$, then $r\in\mathbb Z$.
\end{corollary}

\begin{proof}
If $r\notin\mathbb Z$, then
Lemma~\ref{lem:rational-translate-modular-form} gives a nonzero
arithmetic modular form $F_r$ with $F_r(\rho)=0$, contradicting
Theorem~\ref{thm:finite-level-genericity}.
\end{proof}

\subsection{The two reflection-symmetric zeros}

The zeros on $\operatorname{Re}\tau=0$ and
$\operatorname{Re}\tau=1/2$ are
\cite{Watson1952,ImamogluJermannToth2014}
\begin{align}
 \rho_0&=iy_0,
 &y_0&=0.523521700018\ldots,\label{eq:rho0}\\
 \rho_{1/2}&=\frac12+iy_{1/2},
 &y_{1/2}&=0.130919030397\ldots.\label{eq:rhohalf}
\end{align}
Proposition~\ref{thm:noCM} shows that neither zero is CM.

Watson's cusp-displacement formula \cite{Watson1952} gives exact convergent expressions for their heights:
\begin{equation}\label{eq:Watson-special-heights}
 y_0=
 \left(\frac6\pi+\sum_{r\ge1}(-1)^{r-1}c_r\right)^{-1},
 \quad
 y_{1/2}=\frac14
 \left(\frac6\pi-\sum_{r\ge1}c_r\right)^{-1}.
\end{equation}
The first identity is $\rho_0=\mathcal G(0)$; the second uses $m/n=1/2$, $m'/n=1/2$ and $e^{2\pi ir/2}=(-1)^r$.

\begin{proposition}\label{prop:rational-vertical-degree}
Let
\[
 \rho=a+iy\in\HH,
 \quad a\in\QQ,\quad y\in\RR_{>0},
\]
and suppose that $\rho$ is algebraic. Then
\begin{equation}\label{eq:degree-y2}
 [\QQ(\rho):\QQ]
 =2[\QQ(y^2):\QQ].
\end{equation}
In particular the degree of $\rho$ is even. If in addition $E_2(\rho)=0$, then
\begin{equation}\label{eq:rational-line-degree4}
 [\QQ(\rho):\QQ]\ge4.
\end{equation}
If $R$ is the minimal polynomial of $-y^2=(\rho-a)^2$, then the minimal polynomial of $\rho$ is
\begin{equation}\label{eq:symmetric-minpoly}
 R\bigl((X-a)^2\bigr).
\end{equation}
\end{proposition}

\begin{proof}
Since $a\in\QQ$,
\[
 \QQ(\rho)=\QQ(iy).
\]
The real number $-y^2=(iy)^2$ belongs to this field, and $iy$ is quadratic over $\QQ(y^2)$. It cannot belong to the real field $\QQ(y^2)$, so
\[
 [\QQ(iy):\QQ(y^2)]=2,
\]
which proves \eqref{eq:degree-y2}. If $E_2(\rho)=0$, degree $2$ is impossible by Proposition~\ref{thm:noCM}, hence the degree is at least $4$. If $R$ is the minimal polynomial of $-y^2$, then $R((X-a)^2)$ has degree $2[\QQ(y^2):\QQ]$, vanishes at $\rho$, and has the same degree as the minimal polynomial of $\rho$; hence the two coincide.
\end{proof}

\begin{corollary}\label{cor:symmetric-zero-degree}
If either $\rho_0$ or $\rho_{1/2}$ is algebraic, its degree is even
and at least $4$. In degree $4$, the corresponding $y_0^2$ or
$y_{1/2}^2$ is a real quadratic irrational. The minimal polynomial of
the zero has the form
\[
 R(X^2)
 \quad\text{for }\rho_0,
 \quad
 R\bigl((X-\tfrac12)^2\bigr)
 \quad\text{for }\rho_{1/2},
\]
with $R\in\QQ[X]$ irreducible.
\end{corollary}

The values of $\Delta$ at these zeros satisfy the following consequence
of Nesterenko's theorem.

\begin{corollary}\label{cor:Delta-zero-trans}
For every zero $\rho$ of $E_2$,
\[
 \Delta(\rho)
 =\frac{E_4(\rho)^3-E_6(\rho)^2}{1728}
\]
is transcendental. Moreover
\[
 e^{2\pi i\rho}
 \quad\text{and}\quad
 \Delta(\rho)
\]
are algebraically independent. The same is true of
$e^{2\pi i\rho}$ and $j(\rho)$.
\end{corollary}

\begin{proof}
Corollary~\ref{cor:Nest-zero-scale} gives the algebraic independence of
\[
 q_\rho=e^{2\pi i\rho},\quad E_4(\rho),\quad E_6(\rho).
\]
Put $u=E_4(\rho)$ and $v=E_6(\rho)$.  The field
$\Qbar(u,v)$ has transcendence degree two, and $q_\rho$ is
transcendental over it.  Both
\[
 \Delta(\rho)=\frac{u^3-v^2}{1728},
 \qquad
 j(\rho)=1728\frac{u^3}{u^3-v^2}
\]
are nonconstant rational functions of the algebraically independent pair
\[
 u,\qquad v.
\]
To see this directly for $\Delta(\rho)$, if it were algebraic then
$u^3-v^2$ would be algebraic, giving a nonzero polynomial relation
between $u$ and $v$.  The same argument follows for $j(\rho)$ after
clearing the denominator in its displayed expression.  Thus each of
$\Delta(\rho)$ and $j(\rho)$ is transcendental.  Since
$q_\rho$ is transcendental over the larger field $\Qbar(u,v)$, it is
transcendental over either subfield
$\Qbar(\Delta(\rho))$ or $\Qbar(j(\rho))$.  Each corresponding pair
is therefore algebraically independent.
\end{proof}

For the two symmetric zeros all the quantities in the preceding corollary are real. In particular
\[
 q_0=e^{-2\pi y_0}>0,
 \quad
 q_{1/2}=-e^{-2\pi y_{1/2}}<0
\]
are real transcendental numbers, and each of the triples
\[
 \{q_0,E_4(\rho_0),E_6(\rho_0)\},
 \quad
 \{q_{1/2},E_4(\rho_{1/2}),E_6(\rho_{1/2})\}
\]
 is algebraically independent. Numerically,
\begin{align*}
 j(\rho_0)&\approx 163787.793193,
 &\Delta(\rho_0)&\approx 0.014468341930,\\
 j(\rho_{1/2})&\approx -161723.794971,
 &\Delta(\rho_{1/2})&\approx -59.279808923.
\end{align*}
Since $D\Delta=E_2\Delta$, the values $\Delta(\rho_0)$ and
$\Delta(\rho_{1/2})$ are critical values of $\Delta$.

\subsection{The reciprocal series in the nome}
\label{subsec:watson-nome}

Watson's construction can also be written directly in the nome. Put
\[
 u=e^{2\pi i\tau}
\]
and define the germ at $u=0$
\begin{equation}\label{eq:Watson-psi}
 \psi(u)=\frac1{24}\left(\frac1{E_2(\tau)}-1\right).
\end{equation}
Because $E_2$ is one-periodic, this is a well-defined meromorphic function on the open unit disk. Watson's auxiliary function is precisely this reciprocal form of $E_2$ in his normalization \cite{Watson1952}.

\begin{proposition}\label{prop:psi-natural-boundary}
The poles of $\psi$ accumulate at every point of the unit circle. Consequently $|u|=1$ is a natural boundary for meromorphic continuation of $\psi$. In particular, the germ $\psi$ satisfies no nonzero linear differential equation with coefficients in $\Qbar(u)$; hence it is not a $G$-function. It is not an $E$-function either.
\end{proposition}

\begin{proof}
The poles of $\psi$ are exactly the numbers
\[
 u_{m/n}=e^{2\pi i\rho_{m/n}}
\]
attached to the zeros of $E_2$ in a translation strip. Let $e^{2\pi ix}$ be any point of the unit circle. Choose reduced fractions $m_k/n_k\to x$ with $n_k\to\infty$. Proposition~\ref{prop:Watson-profile} gives
\[
 \rho_{m_k/n_k}
 =\frac{m_k}{n_k}+\frac1{n_k^2}
 \mathcal G\!\left(\frac{m_k'}{n_k}\right).
\]
The cusp-displacement function $\mathcal G$ is bounded on $\RR$, so $\rho_{m_k/n_k}\to x$ from the upper half-plane. Therefore
\[
 u_{m_k/n_k}\longrightarrow e^{2\pi ix}.
\]
Thus every boundary point is an accumulation point of poles, which prevents meromorphic continuation across any point of $|u|=1$.

A solution of a linear differential equation over $\Qbar(u)$ has
finite singularities only at singular points of the equation. The
infinitely many poles of $\psi$ exclude such an equation and hence
exclude the $G$-function property \cite{Andre1989}. They also exclude
the $E$-function property, since an $E$-function is entire.
\end{proof}

The imaginary-axis zero $\rho_0=iy_0$ has maximal height. For every
cusp of denominator $n\ge2$, \eqref{eq:Watson-profile} and
\eqref{eq:G-IJT-bound} give
\[
 \operatorname{Im}\rho_{m/n}
 \leq \frac1{n^2}\left(\frac\pi6+7.77\times10^{-5}\right)
 < y_0.
\]
Hence the radius of convergence of the Taylor series of $\psi$ at $0$ is exactly
\begin{equation}\label{eq:psi-radius}
 R_\psi=e^{-2\pi y_0}.
\end{equation}

\section{Chudnovsky and Nesterenko in elliptic coordinates}

Chudnovsky's theorem has the following form
\cite[Chapter~7, Theorem~3.1, p.~309]{Chudnovsky1984}; see also
\cite[\S5.2]{Waldschmidt2006}.

\begin{theorem}[Chudnovsky]\label{thm:Chudnovsky}
Let $\wp$ be a Weierstrass elliptic function with invariants $g_2,g_3$. Let $\omega\neq0$ be a period and let $\eta=H_\Lambda(\omega)$ be the corresponding quasi-period. Then
\[
 \trdeg_{\QQ}\QQ\left(g_2,g_3,\frac{\omega}{\pi},\frac{\eta}{\omega}\right)\geq2.
\]
Equivalently, at least two of the four numbers
\[
 g_2,\quad g_3,\quad \frac{\omega}{\pi},\quad \frac{\eta}{\omega}
\]
are algebraically independent over $\QQ$.
\end{theorem}

Bertrand's modular consequence is recorded in
\cite[\S5.6]{Waldschmidt2006}.

\begin{proposition}\label{prop:Chud-global}
For every $\tau\in\HH$,
\[
 \trdeg_{\QQ}\QQ\bigl(E_2(\tau),E_4(\tau),E_6(\tau)\bigr)\geq 2.
\]
Equivalently, at least two of $E_2(\tau),E_4(\tau),E_6(\tau)$ are algebraically independent over $\QQ$.
\end{proposition}

\begin{proof}
Apply Theorem~\ref{thm:Chudnovsky} to the lattice
\[
 \Lambda=\pi(\ZZ+\ZZ\tau)
\]
and to the period $\omega=\pi$. By homogeneity and the standard formulas for the Eisenstein invariants,
\[
 g_2(\Lambda)=\frac{4}{3}E_4(\tau),\quad
 g_3(\Lambda)=\frac{8}{27}E_6(\tau),\quad
 \frac{\eta}{\omega}=\frac{1}{3}E_2(\tau),\quad
 \frac{\omega}{\pi}=1.
\]
The assertion follows from Theorem~\ref{thm:Chudnovsky}.
\end{proof}

\begin{corollary}\label{cor:Chud-alginv}
If $g_2,g_3\in\Qbar$ and $\omega\ne0$ is a period with associated quasi-period $\eta$, then
\[
 \frac{\pi}{\omega}\quad\text{and}\quad\frac{\eta}{\omega}
\]
are algebraically independent over $\Qbar$. In particular, $\eta\ne0$.
\end{corollary}

\begin{proof}
When $g_2$ and $g_3$ are algebraic, Theorem~\ref{thm:Chudnovsky} forces $\omega/\pi$ and $\eta/\omega$ to be algebraically independent. Replacing the first nonzero number by its reciprocal preserves algebraic independence.
\end{proof}

\begin{corollary}\label{cor:Chudnovsky-rational}
Assume that $[\eta_1:\eta_2]\in\PP^1(\QQ)$. Then there exists a nonzero period $\omega\in\Lambda$ with $H_\Lambda(\omega)=0$, and
\begin{equation}\label{eq:Chud-bound}
 \trdeg_{\QQ}\QQ\left(g_2,g_3,\frac{\omega}{\pi}\right)\geq2.
\end{equation}
In particular, $g_2$ and $g_3$ cannot both be algebraic. More precisely, if one of $g_2,g_3$ is algebraic, then the other one and $\omega/\pi$ are algebraically independent.
\end{corollary}

\begin{proof}
By Proposition~\ref{prop:rational-zero}, after an integral change of basis there is a primitive period $\omega$ whose associated quasi-period is zero. Apply Theorem~\ref{thm:Chudnovsky}. Since $\eta/\omega=0$ is algebraic, the remaining three quantities must have transcendence degree at least two. The remaining assertions follow.
\end{proof}

\begin{corollary}\label{cor:E4E6-Chudnovsky}
If $\rho\in\HH$ satisfies $E_2(\rho)=0$, then
\[
 E_4(\rho)\quad\text{and}\quad E_6(\rho)
\]
are algebraically independent over $\Qbar$.
\end{corollary}

\begin{proof}
This is the specialization of Proposition~\ref{prop:Chud-global} at $E_2(\rho)=0$. Since one of the three generators is then algebraic, the remaining two must be algebraically independent.
\end{proof}

Corollary~\ref{cor:E4E6-Chudnovsky} also implies the transcendence of $j(\rho)$. Nesterenko's theorem \cite{Nesterenko1996} gives the stronger independence of
\[
 e^{2\pi i\rho},\quad E_4(\rho),\quad E_6(\rho).
\]
For a general lattice with a rational quasi-period relation,
Chudnovsky's theorem gives \eqref{eq:Chud-bound}.

\subsection{Nesterenko's theorem in elliptic coordinates}

\begin{theorem}[Nesterenko]\label{thm:Nesterenko-classical}
For every $\tau\in\HH$, with $q=e^{2\pi i\tau}$,
\[
 \trdeg_{\QQ}\QQ\bigl(q,E_2(\tau),E_4(\tau),E_6(\tau)\bigr)\geq3.
\]
\end{theorem}

For algebraic invariants $g_2,g_3$, the elliptic corollary
\cite[Chapter~3, \S1]{NesterenkoPhilippon2001} asserts the algebraic
independence of
\[
 e^{2\pi i\tau},\quad \frac{\omega}{\pi},\quad
 \frac{\eta}{\pi}.
\]
For arbitrary invariants, one obtains the following estimate.

\begin{proposition}\label{prop:Nest-elliptic}
Let
\[
 \Lambda=\ZZ\omega+\ZZ\omega',\quad
 \tau=\frac{\omega'}{\omega}\in\HH,
\]
let $g_2,g_3$ be the invariants of $\Lambda$, and let
$\eta=H_\Lambda(\omega)$. Put $q=e^{2\pi i\tau}$. Then
\begin{equation}\label{eq:Nest-elliptic}
 \trdeg_{\QQ}
 \QQ\left(g_2,g_3,\frac{\omega}{\pi},
 \frac{\eta}{\pi},q\right)\ge3.
\end{equation}
Equivalently, at least three of the five numbers
\[
 g_2,\quad g_3,\quad \frac{\omega}{\pi},\quad
 \frac{\eta}{\pi},\quad e^{2\pi i\tau}
\]
are algebraically independent over $\QQ$.
\end{proposition}

\begin{proof}
Write $P=E_2$, $Q=E_4$ and $R=E_6$. The standard relations between
Ramanujan's functions \cite{NesterenkoPhilippon2001} and the periods of the Weierstrass function are
\begin{equation}\label{eq:Nest-period-identities}
 P(q)=3\frac{\omega}{\pi}\frac{\eta}{\pi},\quad
 Q(q)=\frac34\left(\frac{\omega}{\pi}\right)^4g_2,
 \quad
 R(q)=\frac{27}{8}\left(\frac{\omega}{\pi}\right)^6g_3.
\end{equation}
Theorem~\ref{thm:Nesterenko-classical} gives
\[
 \trdeg_{\QQ}\QQ\bigl(q,P(q),Q(q),R(q)\bigr)\ge3.
\]
By \eqref{eq:Nest-period-identities}, this field is contained in
\[
 \QQ\left(q,g_2,g_3,\frac{\omega}{\pi},\frac{\eta}{\pi}\right),
\]
which proves \eqref{eq:Nest-elliptic}. The same period identities
are also recalled in \cite[\S5.6]{Waldschmidt2006}.
\end{proof}

If $g_2,g_3\in\Qbar$, Proposition~\ref{prop:Nest-elliptic} reduces to
\[
 \trdeg_{\QQ}\QQ\left(q,\frac{\omega}{\pi},
 \frac{\eta}{\pi}\right)=3,
\]
the usual elliptic corollary of Nesterenko
\cite{Nesterenko1996,NesterenkoPhilippon2001}.

\begin{corollary}\label{cor:Nest-rational}
Assume that $[\eta_1:\eta_2]\in\PP^1(\QQ)$. After an integral change of
basis there are a primitive period $\omega$, a second period $\omega'$, and
\[
 \tau'=\frac{\omega'}{\omega}\in\HH
\]
such that $H_\Lambda(\omega)=0$. For $q'=e^{2\pi i\tau'}$ one has
\begin{equation}\label{eq:Nest-rational}
 \trdeg_{\QQ}\QQ\left(g_2,g_3,\frac{\omega}{\pi},q'\right)\ge3.
\end{equation}
Consequently, if $g_2\in\Qbar$, then
\[
 g_3,\quad \frac{\omega}{\pi},\quad q'
\]
are algebraically independent; the analogous assertion holds with $g_2$
and $g_3$ interchanged. In particular $g_2$ and $g_3$ cannot both be
algebraic.
\end{corollary}

\begin{proof}
Use Proposition~\ref{prop:rational-zero} to choose an integral basis for
which the first quasi-period vanishes, and then apply
Proposition~\ref{prop:Nest-elliptic}. Since $\eta/\pi=0$, the remaining
four generators have transcendence degree at least three. If one of the
invariants is algebraic, the remaining three generators must therefore be
algebraically independent.
\end{proof}

\begin{corollary}\label{cor:Nest-zero-scale}
Let $E_2(\rho)=0$ and let $\omega\in\CC^\times$. For the lattice
\[
 \Lambda_\omega=\omega(\ZZ+\ZZ\rho)
\]
with invariants $g_2(\Lambda_\omega),g_3(\Lambda_\omega)$, one has
\[
 \trdeg_{\QQ}\QQ\left(
 g_2(\Lambda_\omega),g_3(\Lambda_\omega),
 \frac{\omega}{\pi},e^{2\pi i\rho}\right)\ge3.
\]
For $\omega=\pi$ this becomes the algebraic independence of
\[
 e^{2\pi i\rho},\quad E_4(\rho),\quad E_6(\rho).
\]
\end{corollary}

\begin{proof}
The quasi-period corresponding to $\omega$ is zero by homogeneity and
\eqref{eq:eta-E2-general}, so Proposition~\ref{prop:Nest-elliptic} applies.
For $\omega=\pi$ one has
\[
 g_2=\frac43E_4(\rho),\quad g_3=\frac8{27}E_6(\rho),
\]
and the last assertion follows.
\end{proof}

Proposition~\ref{prop:Nest-elliptic} also proves that $j(\rho)$ is
transcendental. Otherwise, a homothetic lattice would have algebraic
$g_2,g_3$ and a vanishing quasi-period, forcing
$\QQ(e^{2\pi i\rho},\omega/\pi)$ to have transcendence degree at
least three.

\section{The Legendre differential module and its Schwarz map}
\label{subsec:watson-G}

Let $\lambda(\tau)$ be the modular lambda function and put
\begin{equation}\label{eq:Legendre-F01}
 F_0(z)={}_2F_1\!\left(\frac12,\frac12;1;z\right),
 \quad
 F_1(z)={}_2F_1\!\left(-\frac12,\frac12;1;z\right).
\end{equation}
These hypergeometric series are classical $G$-functions in the sense of Siegel; see
\cite{Andre1989}. At a point $\tau\in\HH$, the notation
$F_j(\lambda(\tau))$ below refers to the analytic continuation of the Legendre
period branch determined by $\tau$. Thus the first branch is normalized by
\begin{equation}\label{eq:F0-theta}
 F_0(\lambda(\tau))=\vartheta_3(\tau)^2
\end{equation}
where
\[
\vartheta_3(\tau)=\sum_{n\in\ZZ}e^{\pi i n^2\tau}
\]
denotes the third Jacobi theta constant.

For points on the positive imaginary axis this is the ordinary principal branch of
\eqref{eq:Legendre-F01}.

The standard theta-constant identities give, with $z=\lambda(\tau)$,
\begin{align}
 E_4(\tau)&=(1-z+z^2)F_0(z)^4,\label{eq:E4-lambda-F}\\
 E_6(\tau)&=\frac12(1+z)(2-z)(1-2z)F_0(z)^6.\label{eq:E6-lambda-F}
\end{align}
For $E_2$, define
\begin{equation}\label{eq:A-Gfunction}
 A(z)=3F_1(z)-(2-z)F_0(z).
\end{equation}
The contiguous identity
\begin{equation}\label{eq:F1-F0prime}
 F_1(z)=(1-z)F_0(z)+2z(1-z)F_0'(z)
\end{equation}
shows that
\begin{equation}\label{eq:A-F0prime}
 A(z)=(1-2z)F_0(z)+6z(1-z)F_0'(z).
\end{equation}
On the other hand,
\[
 \Delta(\tau)=\frac1{256}z^2(1-z)^2F_0(z)^{12},
 \quad
 \frac1{2\pi i}\frac{dz}{d\tau}
 =\frac12z(1-z)F_0(z)^2.
\]
Taking the logarithmic derivative and using $D\Delta=E_2\Delta$ gives
\begin{equation}\label{eq:E2-Gfunction-factorization}
 E_2(\tau)=F_0(z)A(z),\quad z=\lambda(\tau).
\end{equation}
The theta product formula gives $F_0(\lambda(\tau))\ne0$.
Thus, on the corresponding Legendre branch,
$E_2(\tau)=0$ if and only if $A(\lambda(\tau))=0$.

The function $F_0$ satisfies the hypergeometric equation
\begin{equation}\label{eq:Gauss-Legendre}
 z(1-z)F_0''+(1-2z)F_0'-\frac14F_0=0.
\end{equation}
Substitution of \eqref{eq:A-F0prime} into \eqref{eq:Gauss-Legendre} yields an
explicit second-order equation for $A$:
\begin{equation}\label{eq:A-Goperator}
 z(1-z)(z^2-z+1)A''
 +(1-2z)A'
 +\frac{1+z-z^2}{4}A=0.
\end{equation}
To relate the two equations, write $Y=F_0$ and
\[
 a=1-2z,\qquad b=6z(1-z),\qquad A=aY+bY'.
\]
Using the Gauss equation \eqref{eq:Gauss-Legendre} to eliminate $Y''$ gives
\begin{equation}\label{eq:gauge-matrix}
 \binom{A}{A'}
 =M(z)\binom{Y}{Y'},
 \qquad
 M(z)=
 \begin{pmatrix}
 1-2z & 6z(1-z)\\[2mm]
 -\tfrac12 & 1-2z
 \end{pmatrix}.
\end{equation}
Its determinant is
\begin{equation}\label{eq:gauge-det}
 \det M(z)=z^2-z+1.
\end{equation}
Thus the first-order system attached to \eqref{eq:A-Goperator} is
rationally gauge-equivalent over $\QQ(z)$ to the Legendre system.

The two additional singularities of the scalar equation \eqref{eq:A-Goperator} are apparent. Indeed, the Gauss equation is nonsingular at a root $\alpha$ of $\alpha^2-\alpha+1=0$, so for every local solution $Y$ the function $A=(1-2z)Y+6z(1-z)Y'$ is holomorphic there. On a punctured disk around $\alpha$, the gauge matrix is invertible, so the transforms of two independent Gauss solutions remain independent. Their Wronskian is $\det M$ times the Gauss Wronskian and therefore has a simple zero at $\alpha$. The indicial equation of \eqref{eq:A-Goperator} at $\alpha$ is
\[
 r(r-2)=0,
\]
so the local exponents are $0$ and $2$, consistent with the simple
Wronskian zero. Both solutions extend holomorphically, so the
singularities are apparent. They are the two $\lambda$-values above
$j=0$.

\begin{corollary}[Hypergeometric values at the zeros]
\label{thm:Gfunction-zero-AI}
Let $\rho\in\HH$ satisfy $E_2(\rho)=0$, put
\[
 \alpha=\lambda(\rho),\quad q_\rho=e^{2\pi i\rho},
\]
and choose the Legendre branches of $F_0,F_1$ determined by $\rho$. Then:
\begin{enumerate}[label=(\roman*)]
 \item $\alpha$ and $F_0(\alpha)$ are algebraically independent over $\Qbar$;
 \item $F_0(\alpha)$ and $F_1(\alpha)$ are algebraically independent over $\Qbar$;
 \item the three numbers
 \[
 q_\rho,\quad F_0(\alpha),\quad F_1(\alpha)
 \]
 are algebraically independent over $\Qbar$.
\end{enumerate}
\end{corollary}

\begin{proof}
By Corollary~\ref{cor:E4E6-Chudnovsky}, the two numbers
$E_4(\rho),E_6(\rho)$ are algebraically independent. Equations
\eqref{eq:E4-lambda-F}--\eqref{eq:E6-lambda-F} show that
\[
 \Qbar(E_4(\rho),E_6(\rho))
 \subset
 \Qbar(\alpha,F_0(\alpha)).
\]
The field on the left has transcendence degree $2$, while the field on the
right is generated by two elements. Hence
\[
 \trdeg_{\Qbar}\Qbar(\alpha,F_0(\alpha))=2,
\]
which proves (i).

At a zero, \eqref{eq:E2-Gfunction-factorization} and the nonvanishing of
$F_0(\alpha)$ give $A(\alpha)=0$. Therefore
\begin{equation}\label{eq:alpha-F0F1}
 3F_1(\alpha)=(2-\alpha)F_0(\alpha),
 \quad
 \alpha=2-3\frac{F_1(\alpha)}{F_0(\alpha)}.
\end{equation}
Consequently
\[
 \Qbar(\alpha,F_0(\alpha))
 =\Qbar(F_0(\alpha),F_1(\alpha)),
\]
and (ii) follows from (i).

Corollary~\ref{cor:Nest-zero-scale} gives the algebraic independence of
\[
 q_\rho,\quad E_4(\rho),\quad E_6(\rho).
\]
By \eqref{eq:E4-lambda-F}--\eqref{eq:E6-lambda-F}, their field is contained in
\[
 \Qbar(q_\rho,\alpha,F_0(\alpha)).
\]
Hence this latter field has transcendence degree $3$. Using
\eqref{eq:alpha-F0F1} once more identifies it with
$\Qbar(q_\rho,F_0(\alpha),F_1(\alpha))$, proving (iii).
\end{proof}

These assertions follow from modular algebraic independence at the
transcendental parameter $\alpha=\lambda(\rho)$. For special-value
theorems at algebraic parameters, see \cite{Beukers1993,Andre1989}.

Put
\begin{equation}\label{eq:Legendre-Dop}
 \mathcal D=(1-2z)+6z(1-z)\frac{d}{dz}.
\end{equation}
Thus $A=\mathcal D F_0$. Let $Y_1=F_0$ and let $Y_2$ be a second local solution
of \eqref{eq:Gauss-Legendre}, normalized on the inverse-$\lambda$ branch by
\begin{equation}\label{eq:tau-Schwarz-original}
 \frac{Y_2(z)}{Y_1(z)}=\tau.
\end{equation}
Set
\[
 A_1=\mathcal D Y_1,
 \quad
 A_2=\mathcal D Y_2.
\]
The two functions $A_1,A_2$ solve \eqref{eq:A-Goperator}. The classical
Wronskian formula for the inverse modular lambda function is
\begin{equation}\label{eq:lambda-Wronskian}
 \frac{d\tau}{dz}
 =\frac{1}{\pi i\,z(1-z)F_0(z)^2}.
\end{equation}
Using $Y_2=\tau F_0$, equations \eqref{eq:Legendre-Dop} and
\eqref{eq:lambda-Wronskian} give
\[
 A_2
 =\tau A_1+\frac{6}{\pi iF_0}.
\]
Since $A_1=A=E_2/F_0$, we obtain
\begin{equation}\label{eq:h-Schwarz-Goperator}
 \frac{A_2(z)}{A_1(z)}
 =\tau+\frac{6}{\pi iE_2(\tau)}
 =h(\tau),
 \quad z=\lambda(\tau).
\end{equation}
Hence $h$ is the Schwarz map of \eqref{eq:A-Goperator}.

Watson's transformed equation \cite[\S4]{Watson1952} is, in his notation,
\[
 T-\frac db=\frac{6i}{\pi E_2(T)},
\]
which is equivalent to
\begin{equation}\label{eq:Watson-Schwarz-rational}
 h(T)=\frac db\in\QQ.
\end{equation}
Watson's Lagrange series therefore gives the inverse of this Schwarz
map at rational points.

\begin{corollary}
\label{cor:symmetric-Gfunction}
Let
\[
 \alpha_0=\lambda(\rho_0),
 \quad
 \alpha_{1/2}=\lambda(\rho_{1/2}).
\]
Then both numbers are transcendental. Moreover
\[
 \alpha_0\in(0,1)
\]
is the unique zero on $(0,1)$ of the principal branch of the $G$-function $A$,
while
\begin{equation}\label{eq:lambda-half-circle}
 |\alpha_{1/2}-1|=1.
\end{equation}
Numerically, for the standard modular lambda function,
\[
 \alpha_0=0.961149475381\ldots,
\]
whereas
\[
 \alpha_{1/2}
 =0.000787737945\ldots
 -0.039684447316\ldots\,i.
\]
For each of these two points the corresponding continued Legendre values
$F_0(\alpha)$ and $F_1(\alpha)$ are algebraically independent, and adjoining the
corresponding nome gives a triple of algebraically independent numbers as in
Corollary~\ref{thm:Gfunction-zero-AI}.
\end{corollary}

\begin{proof}
The transcendence and algebraic-independence assertions follow from
Corollary~\ref{thm:Gfunction-zero-AI}. For $\rho_0=iy_0$, the standard lambda
parameter lies in $(0,1)$ and the principal Legendre branch applies; uniqueness
follows from Watson's uniqueness \cite{Watson1952} of the purely imaginary zero. For
$\rho_{1/2}=\frac12+iy_{1/2}$, one has
$-\overline{\rho_{1/2}}=\rho_{1/2}-1$. The real Fourier coefficients of
$\lambda$ and the transformation
\[
 \lambda(\tau-1)=\frac{\lambda(\tau)}{\lambda(\tau)-1}
\]
give
\[
 \overline{\alpha_{1/2}}
 =\frac{\alpha_{1/2}}{\alpha_{1/2}-1}.
\]
Taking absolute values yields \eqref{eq:lambda-half-circle}.
\end{proof}

Write the Taylor expansion at the origin as
\[
 A(z)=\sum_{n\ge0}a_nz^n,
 \quad a_n\in\QQ.
\]
Since $A$ is a $G$-function, its coefficientwise Borel transform
\begin{equation}\label{eq:Efunction-companion}
 \mathcal E_A(w)=\sum_{n\ge0}\frac{a_n}{n!}w^n
\end{equation}
is an $E$-function. The Fourier--Laplace correspondence associates the
minimal $G$-operator of $A$ with an $E$-operator
\cite{Andre2000Gevrey,FischlerRivoal2016}. Special-value theorems for
$E$-functions at algebraic arguments do not apply to the transcendental
parameter $\alpha=\lambda(\rho)$, and Borel transformation does not
preserve the equation $A(\alpha)=0$.

\section{Watson's profile and the real Schwarz locus}
\label{subsec:real-Schwarz-locus}

Regard $h$ as a meromorphic map $\HH\to\PP^1(\CC)$, with
$h(\rho)=\infty$ when $E_2(\rho)=0$. Let $\mathcal F$ be the
standard fundamental domain and put
\[
 \mathcal A
 =
 \{\tau\in\mathcal F:\ h(\tau)\in\RR\}.
\]
Wood and Young \cite[Theorems~5.1 and 6.2]{WoodYoung2014} locate
this arc near $\operatorname{Im}\tau=6/\pi$ and prove uniqueness of
the inverse real branch. Imamo\u{g}lu, Jermann and T\'oth
\cite[Theorem~1.2 and \S2]{ImamogluJermannToth2014} prove that
$\mathcal A$ is a graph over $[-1/2,1/2]$ and that
\[
 h:\mathcal A\longrightarrow[-1/2,1/2]
\]
is bijective, with nonzero derivative. Holomorphic inversion makes
its inverse real analytic. Write
\[
 F_h:[-1/2,1/2]\longrightarrow\mathcal A,
 \quad
 h(F_h(x))=x,
\]
and extend it by $F_h(x+1)=F_h(x)+1$. Thus
\[
 H(x)=F_h(x)-x
\]
is one-periodic. In their notation, $F_h$ is the inverse function $F$. They obtain a Fourier expansion for
$H$, with
\[
 c_0=\frac{6i}{\pi},
 \quad
 c_1=\frac{144e^{-12}}{\pi}\,i
 =0.000281630\ldots\,i,
\]
and the uniform annular estimate
\begin{equation}\label{eq:IJT-H-annulus}
 0.00027
 <
 \left|H(x)-\frac{6i}{\pi}\right|
 <
 0.00029
 \quad(x\in\RR).
\end{equation}
The following identity relates this inverse to Watson's series.

\begin{theorem}\label{thm:equivariant-Watson-dictionary}
For every real $x$,
\begin{equation}\label{eq:Fh-Watson}
 F_h(x)
 =
 x+iP_W(e^{2\pi ix}),
 \quad
 H(x)=iP_W(e^{2\pi ix}),
\end{equation}
and therefore
\begin{equation}\label{eq:G-reciprocal-H}
 \mathcal G(x)=-\frac1{H(x)}.
\end{equation}
More precisely, let
\[
 \gamma=\begin{pmatrix}a&b\\c&d\end{pmatrix}\in\mathrm{SL}_2(\ZZ),
 \quad c>0,
\]
and put
\[
 x=-\frac dc\in\RR.
\]
Here $F_h$ is the translation-equivariant extension to all of $\RR$.
Then
\begin{equation}\label{eq:zero-from-real-locus}
 \rho_{a/c}
 =
 \gamma F_h(x)
 =
 \frac ac-\frac1{c^2H(x)}
 =
 \frac ac+\frac1{c^2}\mathcal G(x).
\end{equation}
If $x_0=x+k$ is the centered representative, with $k\in\ZZ$, put
$T_1=\left(\begin{smallmatrix}1&1\\0&1\end{smallmatrix}\right)$ and
$\gamma_0=\gamma T_1^{-k}$. Then
\[
 \gamma_0F_h(x_0)=\gamma F_h(x)=\rho_{a/c}.
\]
Thus centering $x$ requires replacing $\gamma$ by $\gamma_0$.
The functions $H$ and $\mathcal G$ are one-periodic, whereas
$F_h(x+k)=F_h(x)+k$.
\end{theorem}

\begin{proof}
For rational $x$, Watson's transformed equation is precisely
$h(T)=x$, while \eqref{eq:Watson-T} gives
\[
 T=x+iP_W(e^{2\pi ix}).
\]
The uniqueness theorem of Imamo\u{g}lu, Jermann and T\'oth on the real arc
therefore identifies $T$ with $F_h(x)$. Since the rational points are dense
and both sides of \eqref{eq:Fh-Watson} are real-analytic, the identity holds
for every real $x$. Equation \eqref{eq:G-reciprocal-H} follows
from the definition \eqref{eq:Watson-profile-function}.

Now put $x=-d/c$ exactly and $z=F_h(x)=x+H(x)$. Since
$ad-bc=1$,
\[
 ax+b=-\frac{ad}{c}+b=-\frac1c,
 \quad
 cx+d=0.
\]
Hence
\[
 \gamma z
 =
 \frac{-1/c+aH(x)}{cH(x)}
 =
 \frac ac-\frac1{c^2H(x)}.
\]
Equivariance gives
\[
 h(\gamma z)=\gamma h(z)=\gamma(-d/c)=\infty,
\]
so $\gamma z$ is a zero of $E_2$. Substitution of
\eqref{eq:G-reciprocal-H} gives the last expression in
\eqref{eq:zero-from-real-locus}. Also,
$ad\equiv1\pmod c$, so $-d$ is the inverse residue characterized by
\eqref{eq:Watson-inverse}. Finally, $F_h(x+k)=T_1^kF_h(x)$ gives
$\gamma T_1^{-k}F_h(x+k)=\gamma F_h(x)$, which proves the assertion
about centered representatives.
\end{proof}

Equivariance also gives the following description of the real locus;
compare \cite[Sections~4--6]{WoodYoung2014}.

\begin{proposition}\label{prop:modular-real-Schwarz-locus}
Let
\[
 \mathcal C=h^{-1}\!\bigl(\PP^1(\RR)\bigr).
\]
Then
\begin{equation}\label{eq:real-locus-union}
 \mathcal C
 =
 \bigcup_{\gamma\in\mathrm{SL}_2(\ZZ)}\gamma\mathcal A.
\end{equation}
Moreover
\begin{equation}\label{eq:rational-Schwarz-locus}
 h^{-1}\!\bigl(\PP^1(\QQ)\bigr)
 =
 \mathrm{SL}_2(\ZZ)\cdot Z(E_2),
\end{equation}
where $Z(E_2)$ is the zero set of $E_2$. Consequently
\begin{equation}\label{eq:dense-zero-saturation}
 \overline{\mathrm{SL}_2(\ZZ)\cdot Z(E_2)}
 =
 \mathcal C.
\end{equation}
\end{proposition}

\begin{proof}
Equivariance implies that $\mathcal C$ is $\mathrm{SL}_2(\ZZ)$-invariant.
Every point of $\HH$ can be moved to $\mathcal F$. On $\mathcal F$ the
function $E_2$ has no zero, and the analysis of
\cite[\S2]{ImamogluJermannToth2014} shows that the points at which $h$ is real
are exactly $\mathcal A$. This proves \eqref{eq:real-locus-union}.

If $\rho$ is a zero, then $h(\rho)=\infty$, and hence
\[
 h(\gamma\rho)=\gamma\infty\in\PP^1(\QQ).
\]
Conversely, if $h(\tau)=r\in\PP^1(\QQ)$, choose
$\gamma\in\mathrm{SL}_2(\ZZ)$ with $\gamma r=\infty$. Then
$h(\gamma\tau)=\infty$, so $\gamma\tau$ is a zero. This proves
\eqref{eq:rational-Schwarz-locus}. On each arc
$\gamma\mathcal A$ the projective map $h$ is a homeomorphism onto a real
projective interval, and $\PP^1(\QQ)$ is dense in $\PP^1(\RR)$.
Equation \eqref{eq:dense-zero-saturation} follows.
\end{proof}

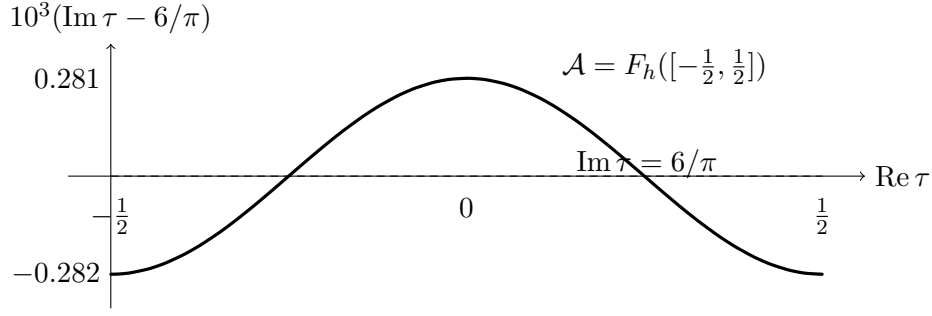
\begin{figure}[t]
\centering
\begin{tikzpicture}[x=9.4cm,y=4.6cm]
 \draw[->] (-0.56,0) -- (0.56,0) node[right] {$\operatorname{Re}\tau$};
 \draw[->] (-0.5,-0.38) -- (-0.5,0.38)
 node[above] {$10^3(\operatorname{Im}\tau-6/\pi)$};
 \draw[densely dashed] (-0.5,0) -- (0.5,0);
 \draw[very thick, smooth] plot coordinates {(-0.5000000,-0.2820827) (-0.4899823,-0.2815234) (-0.4799646,-0.2798477) (-0.4699471,-0.2770623) (-0.4499127,-0.2682120) (-0.4398960,-0.2621826) (-0.4298797,-0.2551145) (-0.4198639,-0.2470363) (-0.3998340,-0.2279826) (-0.3898200,-0.2170837) (-0.3798068,-0.2053271) (-0.3697942,-0.1927598) (-0.3497717,-0.1653973) (-0.3397618,-0.1507115) (-0.3297528,-0.1354332) (-0.3197448,-0.1196233) (-0.2997319,-0.0866622) (-0.2897270,-0.0696420) (-0.2797232,-0.0523517) (-0.2697205,-0.0348598) (-0.2497184,0.0004516) (-0.2397190,0.0181316) (-0.2297207,0.0357346) (-0.2197235,0.0531915) (-0.1997324,0.0873929) (-0.1897385,0.1040031) (-0.1797455,0.1201990) (-0.1697536,0.1359172) (-0.1497726,0.1656765) (-0.1397834,0.1796014) (-0.1297952,0.1928165) (-0.1198077,0.2052704) (-0.0998349,0.2277035) (-0.0898495,0.2375956) (-0.0798647,0.2465523) (-0.0698804,0.2545388) (-0.0499132,0.2674813) (-0.0399302,0.2723870) (-0.0299474,0.2762225) (-0.0199648,0.2789728) (0.0000000,0.2811794) (0.0099824,0.2806272) (0.0199648,0.2789728) (0.0299474,0.2762225) (0.0499132,0.2674813) (0.0598966,0.2615241) (0.0698804,0.2545388) (0.0798647,0.2465523) (0.0998349,0.2277035) (0.1098209,0.2169145) (0.1198077,0.2052704) (0.1297952,0.1928165) (0.1497726,0.1656765) (0.1597626,0.1510961) (0.1697536,0.1359172) (0.1797455,0.1201990) (0.1997324,0.0873929) (0.2097274,0.0704335) (0.2197235,0.0531915) (0.2297207,0.0357346) (0.2497184,0.0004516) (0.2597189,-0.0172354) (0.2697205,-0.0348598) (0.2797232,-0.0523517) (0.2997319,-0.0866622) (0.3097378,-0.1033447) (0.3197448,-0.1196233) (0.3297528,-0.1354332) (0.3497717,-0.1653973) (0.3597826,-0.1794321) (0.3697942,-0.1927598) (0.3798068,-0.2053271) (0.3998340,-0.2279826) (0.4098487,-0.2379802) (0.4198639,-0.2470363) (0.4298797,-0.2551145) (0.4499127,-0.2682120) (0.4599297,-0.2731786) (0.4699471,-0.2770623) (0.4799646,-0.2798477) (0.5000000,-0.2820827)};
 \fill (0,0.2811794) circle (0.55pt);
 \fill (-0.5,-0.2820827) circle (0.55pt);
 \fill (0.5,-0.2820827) circle (0.55pt);
 \node[below] at (-0.5,-0.035) {$-\tfrac12$};
 \node[below] at (0,-0.035) {$0$};
 \node[below] at (0.5,-0.035) {$\tfrac12$};
 \node[left] at (-0.5,0.281) {$0.281$};
 \node[left] at (-0.5,-0.282) {$-0.282$};
 \node[anchor=west] at (0.12,0.315) {$\mathcal A=F_h([-\tfrac12,\tfrac12])$};
 \node[anchor=west] at (0.14,0.035) {$\operatorname{Im}\tau=6/\pi$};
\end{tikzpicture}
\caption{The real inverse arc $\mathcal A=h^{-1}(\RR)$ in the standard fundamental strip. The vertical displacement from $\operatorname{Im}\tau=6/\pi$ is magnified by a factor $10^3$.}
\label{fig:real-Schwarz-arc}
\end{figure}

\begin{figure}[t]
\centering
\begin{tikzpicture}[x=10.0cm,y=10.0cm]
 \draw[->] (-0.57,0) -- (0.57,0) node[right] {$\operatorname{Re}\tau$};
 \draw[->] (0,0) -- (0,0.57) node[above] {$\operatorname{Im}\tau$};
 \draw[densely dashed,black!50] (0.0000000,0.2617994) circle [radius=0.2617994];
 \draw[densely dashed,black!50] (0.5000000,0.0654498) circle [radius=0.0654498];
 \draw[densely dashed,black!50] (-0.5000000,0.0654498) circle [radius=0.0654498];
 \draw[densely dashed,black!50] (0.3333333,0.0290888) circle [radius=0.0290888];
 \draw[densely dashed,black!50] (-0.3333333,0.0290888) circle [radius=0.0290888];
 \draw[densely dashed,black!50] (0.2500000,0.0163625) circle [radius=0.0163625];
 \draw[densely dashed,black!50] (-0.2500000,0.0163625) circle [radius=0.0163625];
 \draw[densely dashed,black!50] (0.2000000,0.0104720) circle [radius=0.0104720];
 \draw[densely dashed,black!50] (-0.2000000,0.0104720) circle [radius=0.0104720];
 \draw[densely dashed,black!50] (0.4000000,0.0104720) circle [radius=0.0104720];
 \draw[densely dashed,black!50] (-0.4000000,0.0104720) circle [radius=0.0104720];
 \draw[thick,smooth] plot coordinates {(0.1283208,0.4900769) (0.1238024,0.4925583) (0.1192129,0.4949609) (0.1145542,0.4972819) (0.1098285,0.4995187) (0.1050382,0.5016686) (0.1001856,0.5037294) (0.0952733,0.5056985) (0.0903039,0.5075739) (0.0852802,0.5093532) (0.0802050,0.5110346) (0.0750812,0.5126162) (0.0699118,0.5140962) (0.0646998,0.5154730) (0.0594483,0.5167450) (0.0541605,0.5179110) (0.0488395,0.5189696) (0.0434886,0.5199198) (0.0381111,0.5207604) (0.0327102,0.5214908) (0.0272892,0.5221100) (0.0218515,0.5226175) (0.0164004,0.5230128) (0.0109392,0.5232954) (0.0054713,0.5234651) (0.0000000,0.5235217) (-0.0054713,0.5234651) (-0.0109392,0.5232954) (-0.0164004,0.5230128) (-0.0218515,0.5226175) (-0.0272892,0.5221100) (-0.0327102,0.5214908) (-0.0381111,0.5207604) (-0.0434886,0.5199198) (-0.0488395,0.5189696) (-0.0541605,0.5179110) (-0.0594483,0.5167450) (-0.0646998,0.5154730) (-0.0699118,0.5140962) (-0.0750812,0.5126162) (-0.0802050,0.5110346) (-0.0852802,0.5093532) (-0.0903039,0.5075739) (-0.0952733,0.5056985) (-0.1001856,0.5037294) (-0.1050382,0.5016686) (-0.1098285,0.4995187) (-0.1145542,0.4972819) (-0.1192129,0.4949609) (-0.1238024,0.4925583) (-0.1283208,0.4900769)};
 \draw[thick,smooth] plot coordinates {(0.5538041,0.1027431) (0.5531800,0.1036278) (0.5525357,0.1045092) (0.5518710,0.1053870) (0.5511858,0.1062603) (0.5504800,0.1071288) (0.5497534,0.1079916) (0.5490061,0.1088483) (0.5482379,0.1096982) (0.5474490,0.1105407) (0.5466392,0.1113751) (0.5458086,0.1122008) (0.5449573,0.1130172) (0.5440853,0.1138237) (0.5431927,0.1146195) (0.5422797,0.1154040) (0.5413463,0.1161767) (0.5403927,0.1169367) (0.5394191,0.1176836) (0.5384257,0.1184167) (0.5374127,0.1191352) (0.5363803,0.1198386) (0.5353288,0.1205262) (0.5342584,0.1211973) (0.5331696,0.1218514) (0.5320625,0.1224877) (0.5309376,0.1231057) (0.5297952,0.1237046) (0.5286358,0.1242839) (0.5274597,0.1248429) (0.5262674,0.1253810) (0.5250595,0.1258976) (0.5238364,0.1263920) (0.5225987,0.1268637) (0.5213470,0.1273121) (0.5200818,0.1277366) (0.5188038,0.1281366) (0.5175138,0.1285116) (0.5162123,0.1288612) (0.5149001,0.1291847) (0.5135781,0.1294818) (0.5122469,0.1297519) (0.5109075,0.1299948) (0.5095606,0.1302100) (0.5082072,0.1303972) (0.5068481,0.1305561) (0.5054842,0.1306865) (0.5041165,0.1307881) (0.5027459,0.1308608) (0.5013735,0.1309045) (0.5000000,0.1309190)};
 \draw[thick,smooth] plot coordinates {(-0.4461959,0.1027431) (-0.4468200,0.1036278) (-0.4474643,0.1045092) (-0.4481290,0.1053870) (-0.4488142,0.1062603) (-0.4495200,0.1071288) (-0.4502466,0.1079916) (-0.4509939,0.1088483) (-0.4517621,0.1096982) (-0.4525510,0.1105407) (-0.4533608,0.1113751) (-0.4541914,0.1122008) (-0.4550427,0.1130172) (-0.4559147,0.1138237) (-0.4568073,0.1146195) (-0.4577203,0.1154040) (-0.4586537,0.1161767) (-0.4596073,0.1169367) (-0.4605809,0.1176836) (-0.4615743,0.1184167) (-0.4625873,0.1191352) (-0.4636197,0.1198386) (-0.4646712,0.1205262) (-0.4657416,0.1211973) (-0.4668304,0.1218514) (-0.4679375,0.1224877) (-0.4690624,0.1231057) (-0.4702048,0.1237046) (-0.4713642,0.1242839) (-0.4725403,0.1248429) (-0.4737326,0.1253810) (-0.4749405,0.1258976) (-0.4761636,0.1263920) (-0.4774013,0.1268637) (-0.4786530,0.1273121) (-0.4799182,0.1277366) (-0.4811962,0.1281366) (-0.4824862,0.1285116) (-0.4837877,0.1288612) (-0.4850999,0.1291847) (-0.4864219,0.1294818) (-0.4877531,0.1297519) (-0.4890925,0.1299948) (-0.4904394,0.1302100) (-0.4917928,0.1303972) (-0.4931519,0.1305561) (-0.4945158,0.1306865) (-0.4958835,0.1307881) (-0.4972541,0.1308608) (-0.4986265,0.1309045) (-0.5000000,0.1309190)};
 \draw[thick,smooth] plot coordinates {(0.3383734,0.0577463) (0.3377751,0.0578451) (0.3371739,0.0579314) (0.3365702,0.0580050) (0.3359645,0.0580659) (0.3353571,0.0581141) (0.3347485,0.0581495) (0.3341390,0.0581719) (0.3335292,0.0581816) (0.3329194,0.0581783) (0.3323100,0.0581623) (0.3317015,0.0581334) (0.3310942,0.0580919) (0.3304886,0.0580377) (0.3298850,0.0579709) (0.3292838,0.0578917) (0.3286854,0.0578002) (0.3280902,0.0576966) (0.3274985,0.0575809) (0.3269106,0.0574534) (0.3263270,0.0573142) (0.3257478,0.0571635) (0.3251735,0.0570016) (0.3246043,0.0568285) (0.3240405,0.0566446) (0.3234824,0.0564500) (0.3229303,0.0562450) (0.3223844,0.0560298) (0.3218450,0.0558046) (0.3213124,0.0555697) (0.3207867,0.0553254) (0.3202682,0.0550718) (0.3197571,0.0548093) (0.3192535,0.0545381) (0.3187578,0.0542586) (0.3182701,0.0539708) (0.3177905,0.0536753) (0.3173192,0.0533721) (0.3168564,0.0530617) (0.3164022,0.0527443) (0.3159567,0.0524202) (0.3155201,0.0520897) (0.3150925,0.0517531) (0.3146739,0.0514107) (0.3142645,0.0510628) (0.3138643,0.0507098) (0.3134733,0.0503519) (0.3130917,0.0499894) (0.3127195,0.0496228) (0.3123567,0.0492521) (0.3120032,0.0488778)};
 \draw[thick,smooth] plot coordinates {(-0.3120032,0.0488778) (-0.3123567,0.0492521) (-0.3127195,0.0496228) (-0.3130917,0.0499894) (-0.3134733,0.0503519) (-0.3138643,0.0507098) (-0.3142645,0.0510628) (-0.3146739,0.0514107) (-0.3150925,0.0517531) (-0.3155201,0.0520897) (-0.3159567,0.0524202) (-0.3164022,0.0527443) (-0.3168564,0.0530617) (-0.3173192,0.0533721) (-0.3177905,0.0536753) (-0.3182701,0.0539708) (-0.3187578,0.0542586) (-0.3192535,0.0545381) (-0.3197571,0.0548093) (-0.3202682,0.0550718) (-0.3207867,0.0553254) (-0.3213124,0.0555697) (-0.3218450,0.0558046) (-0.3223844,0.0560298) (-0.3229303,0.0562450) (-0.3234824,0.0564500) (-0.3240405,0.0566446) (-0.3246043,0.0568285) (-0.3251735,0.0570016) (-0.3257478,0.0571635) (-0.3263270,0.0573142) (-0.3269106,0.0574534) (-0.3274985,0.0575809) (-0.3280902,0.0576966) (-0.3286854,0.0578002) (-0.3292838,0.0578917) (-0.3298850,0.0579709) (-0.3304886,0.0580377) (-0.3310942,0.0580919) (-0.3317015,0.0581334) (-0.3323100,0.0581623) (-0.3329194,0.0581783) (-0.3335292,0.0581816) (-0.3341390,0.0581719) (-0.3347485,0.0581495) (-0.3353571,0.0581141) (-0.3359645,0.0580659) (-0.3365702,0.0580050) (-0.3371739,0.0579314) (-0.3377751,0.0578451) (-0.3383734,0.0577463)};
 \draw[thick,smooth] plot coordinates {(0.2542127,0.0321782) (0.2538852,0.0322618) (0.2535552,0.0323388) (0.2532228,0.0324089) (0.2528885,0.0324722) (0.2525523,0.0325286) (0.2522145,0.0325779) (0.2518753,0.0326202) (0.2515349,0.0326554) (0.2511936,0.0326834) (0.2508516,0.0327042) (0.2505092,0.0327179) (0.2501665,0.0327244) (0.2498238,0.0327237) (0.2494814,0.0327158) (0.2491394,0.0327008) (0.2487980,0.0326787) (0.2484575,0.0326495) (0.2481181,0.0326133) (0.2477800,0.0325701) (0.2474434,0.0325200) (0.2471085,0.0324631) (0.2467755,0.0323995) (0.2464445,0.0323292) (0.2461159,0.0322524) (0.2457897,0.0321691) (0.2454661,0.0320794) (0.2451454,0.0319834) (0.2448276,0.0318813) (0.2445130,0.0317731) (0.2442016,0.0316590) (0.2438937,0.0315392) (0.2435894,0.0314136) (0.2432888,0.0312826) (0.2429921,0.0311461) (0.2426994,0.0310044) (0.2424109,0.0308576) (0.2421266,0.0307058) (0.2418467,0.0305492) (0.2415713,0.0303880) (0.2413005,0.0302222) (0.2410343,0.0300522) (0.2407730,0.0298780) (0.2405165,0.0296998) (0.2402650,0.0295178) (0.2400185,0.0293322) (0.2397770,0.0291431) (0.2395407,0.0289507) (0.2393096,0.0287552) (0.2390837,0.0285568) (0.2388631,0.0283557)};
 \draw[thick,smooth] plot coordinates {(-0.2388631,0.0283557) (-0.2390837,0.0285568) (-0.2393096,0.0287552) (-0.2395407,0.0289507) (-0.2397770,0.0291431) (-0.2400185,0.0293322) (-0.2402650,0.0295178) (-0.2405165,0.0296998) (-0.2407730,0.0298780) (-0.2410343,0.0300522) (-0.2413005,0.0302222) (-0.2415713,0.0303880) (-0.2418467,0.0305492) (-0.2421266,0.0307058) (-0.2424109,0.0308576) (-0.2426994,0.0310044) (-0.2429921,0.0311461) (-0.2432888,0.0312826) (-0.2435894,0.0314136) (-0.2438937,0.0315392) (-0.2442016,0.0316590) (-0.2445130,0.0317731) (-0.2448276,0.0318813) (-0.2451454,0.0319834) (-0.2454661,0.0320794) (-0.2457897,0.0321691) (-0.2461159,0.0322524) (-0.2464445,0.0323292) (-0.2467755,0.0323995) (-0.2471085,0.0324631) (-0.2474434,0.0325200) (-0.2477800,0.0325701) (-0.2481181,0.0326133) (-0.2484575,0.0326495) (-0.2487980,0.0326787) (-0.2491394,0.0327008) (-0.2494814,0.0327158) (-0.2498238,0.0327237) (-0.2501665,0.0327244) (-0.2505092,0.0327179) (-0.2508516,0.0327042) (-0.2511936,0.0326834) (-0.2515349,0.0326554) (-0.2518753,0.0326202) (-0.2522145,0.0325779) (-0.2525523,0.0325286) (-0.2528885,0.0324722) (-0.2532228,0.0324089) (-0.2535552,0.0323388) (-0.2538852,0.0322618) (-0.2542127,0.0321782)};
 \draw[thick,smooth] plot coordinates {(0.2032116,0.0204425) (0.2030064,0.0205062) (0.2027994,0.0205659) (0.2025906,0.0206214) (0.2023802,0.0206726) (0.2021682,0.0207196) (0.2019549,0.0207621) (0.2017404,0.0208003) (0.2015247,0.0208340) (0.2013082,0.0208632) (0.2010908,0.0208879) (0.2008728,0.0209081) (0.2006543,0.0209237) (0.2004354,0.0209347) (0.2002163,0.0209411) (0.1999971,0.0209430) (0.1997779,0.0209403) (0.1995590,0.0209330) (0.1993404,0.0209212) (0.1991223,0.0209049) (0.1989048,0.0208840) (0.1986881,0.0208587) (0.1984722,0.0208290) (0.1982574,0.0207949) (0.1980436,0.0207565) (0.1978312,0.0207138) (0.1976201,0.0206668) (0.1974105,0.0206157) (0.1972025,0.0205604) (0.1969962,0.0205011) (0.1967918,0.0204378) (0.1965893,0.0203706) (0.1963888,0.0202996) (0.1961905,0.0202248) (0.1959945,0.0201462) (0.1958008,0.0200641) (0.1956095,0.0199785) (0.1954208,0.0198895) (0.1952346,0.0197971) (0.1950512,0.0197014) (0.1948706,0.0196027) (0.1946927,0.0195009) (0.1945179,0.0193961) (0.1943460,0.0192886) (0.1941771,0.0191783) (0.1940114,0.0190654) (0.1938488,0.0189500) (0.1936894,0.0188322) (0.1935333,0.0187121) (0.1933805,0.0185899) (0.1932310,0.0184657)};
 \draw[thick,smooth] plot coordinates {(-0.1932310,0.0184657) (-0.1933805,0.0185899) (-0.1935333,0.0187121) (-0.1936894,0.0188322) (-0.1938488,0.0189500) (-0.1940114,0.0190654) (-0.1941771,0.0191783) (-0.1943460,0.0192886) (-0.1945179,0.0193961) (-0.1946927,0.0195009) (-0.1948706,0.0196027) (-0.1950512,0.0197014) (-0.1952346,0.0197971) (-0.1954208,0.0198895) (-0.1956095,0.0199785) (-0.1958008,0.0200641) (-0.1959945,0.0201462) (-0.1961905,0.0202248) (-0.1963888,0.0202996) (-0.1965893,0.0203706) (-0.1967918,0.0204378) (-0.1969962,0.0205011) (-0.1972025,0.0205604) (-0.1974105,0.0206157) (-0.1976201,0.0206668) (-0.1978312,0.0207138) (-0.1980436,0.0207565) (-0.1982574,0.0207949) (-0.1984722,0.0208290) (-0.1986881,0.0208587) (-0.1989048,0.0208840) (-0.1991223,0.0209049) (-0.1993404,0.0209212) (-0.1995590,0.0209330) (-0.1997779,0.0209403) (-0.1999971,0.0209430) (-0.2002163,0.0209411) (-0.2004354,0.0209347) (-0.2006543,0.0209237) (-0.2008728,0.0209081) (-0.2010908,0.0208879) (-0.2013082,0.0208632) (-0.2015247,0.0208340) (-0.2017404,0.0208003) (-0.2019549,0.0207621) (-0.2021682,0.0207196) (-0.2023802,0.0206726) (-0.2025906,0.0206214) (-0.2027994,0.0205659) (-0.2030064,0.0205062) (-0.2032116,0.0204425)};
 \draw[thick,smooth] plot coordinates {(0.4080781,0.0171398) (0.4079620,0.0172780) (0.4078425,0.0174153) (0.4077197,0.0175515) (0.4075934,0.0176865) (0.4074639,0.0178203) (0.4073309,0.0179526) (0.4071947,0.0180835) (0.4070550,0.0182127) (0.4069121,0.0183402) (0.4067658,0.0184659) (0.4066163,0.0185897) (0.4064635,0.0187114) (0.4063075,0.0188310) (0.4061483,0.0189484) (0.4059860,0.0190634) (0.4058206,0.0191759) (0.4056521,0.0192859) (0.4054806,0.0193932) (0.4053061,0.0194977) (0.4051288,0.0195993) (0.4049486,0.0196980) (0.4047657,0.0197937) (0.4045801,0.0198861) (0.4043918,0.0199753) (0.4042010,0.0200611) (0.4040077,0.0201435) (0.4038120,0.0202224) (0.4036141,0.0202976) (0.4034139,0.0203691) (0.4032116,0.0204368) (0.4030073,0.0205006) (0.4028011,0.0205604) (0.4025931,0.0206161) (0.4023834,0.0206678) (0.4021721,0.0207152) (0.4019594,0.0207584) (0.4017454,0.0207973) (0.4015301,0.0208317) (0.4013138,0.0208618) (0.4010966,0.0208873) (0.4008787,0.0209083) (0.4006600,0.0209248) (0.4004409,0.0209366) (0.4002215,0.0209439) (0.4000018,0.0209465) (0.3997821,0.0209444) (0.3995626,0.0209377) (0.3993433,0.0209263) (0.3991244,0.0209104) (0.3989061,0.0208898)};
 \draw[thick,smooth] plot coordinates {(-0.3989061,0.0208898) (-0.3991244,0.0209104) (-0.3993433,0.0209263) (-0.3995626,0.0209377) (-0.3997821,0.0209444) (-0.4000018,0.0209465) (-0.4002215,0.0209439) (-0.4004409,0.0209366) (-0.4006600,0.0209248) (-0.4008787,0.0209083) (-0.4010966,0.0208873) (-0.4013138,0.0208618) (-0.4015301,0.0208317) (-0.4017454,0.0207973) (-0.4019594,0.0207584) (-0.4021721,0.0207152) (-0.4023834,0.0206678) (-0.4025931,0.0206161) (-0.4028011,0.0205604) (-0.4030073,0.0205006) (-0.4032116,0.0204368) (-0.4034139,0.0203691) (-0.4036141,0.0202976) (-0.4038120,0.0202224) (-0.4040077,0.0201435) (-0.4042010,0.0200611) (-0.4043918,0.0199753) (-0.4045801,0.0198861) (-0.4047657,0.0197937) (-0.4049486,0.0196980) (-0.4051288,0.0195993) (-0.4053061,0.0194977) (-0.4054806,0.0193932) (-0.4056521,0.0192859) (-0.4058206,0.0191759) (-0.4059860,0.0190634) (-0.4061483,0.0189484) (-0.4063075,0.0188310) (-0.4064635,0.0187114) (-0.4066163,0.0185897) (-0.4067658,0.0184659) (-0.4069121,0.0183402) (-0.4070550,0.0182127) (-0.4071947,0.0180835) (-0.4073309,0.0179526) (-0.4074639,0.0178203) (-0.4075934,0.0176865) (-0.4077197,0.0175515) (-0.4078425,0.0174153) (-0.4079620,0.0172780) (-0.4080781,0.0171398)};
 \fill (0.0000000,0.5235217) circle (0.48pt);
 \fill (0.5000000,0.1309190) circle (0.48pt);
 \fill (-0.5000000,0.1309190) circle (0.48pt);
 \fill (0.3333259,0.0581819) circle (0.48pt);
 \fill (-0.3333259,0.0581819) circle (0.48pt);
 \fill (0.2499952,0.0327249) circle (0.48pt);
 \fill (-0.2499952,0.0327249) circle (0.48pt);
 \fill (0.1999971,0.0209430) circle (0.48pt);
 \fill (-0.1999971,0.0209430) circle (0.48pt);
 \fill (0.4000018,0.0209465) circle (0.48pt);
 \fill (-0.4000018,0.0209465) circle (0.48pt);
 \node[anchor=west] at (0.0250000,0.5355217) {$\rho_0$};
 \node[anchor=west] at (0.5120000,0.1429190) {$\rho_{1/2}$};
 \node[anchor=west] at (-0.4533259,0.0721819) {$\rho_{-1/3}$};
 \node[anchor=west] at (0.3453259,0.0721819) {$\rho_{1/3}$};
 \draw (-0.5,0.005) -- (-0.5,-0.005) node[below] {\scriptsize $-\tfrac12$};
 \draw (-0.3333333,0.004) -- (-0.3333333,-0.004) node[below] {\scriptsize $-\tfrac13$};
 \draw (0.3333333,0.004) -- (0.3333333,-0.004) node[below] {\scriptsize $\tfrac13$};
 \draw (0.5,0.005) -- (0.5,-0.005) node[below] {\scriptsize $\tfrac12$};
 \node[anchor=west] at (0.17,0.555) {\small solid: $\gamma\mathcal A$};
 \node[anchor=west] at (0.17,0.532) {\small dashed: reference horocycles};
\end{tikzpicture}
\caption{Translates $\gamma\mathcal A$ for cusps of denominator at
most $5$ (solid), reference horocycles (dashed), and zeros of $E_2$
(points). The circle based at $a/c$ has highest point
$a/c+\pi i/(6c^2)$.}
\label{fig:deformed-Ford}
\end{figure}
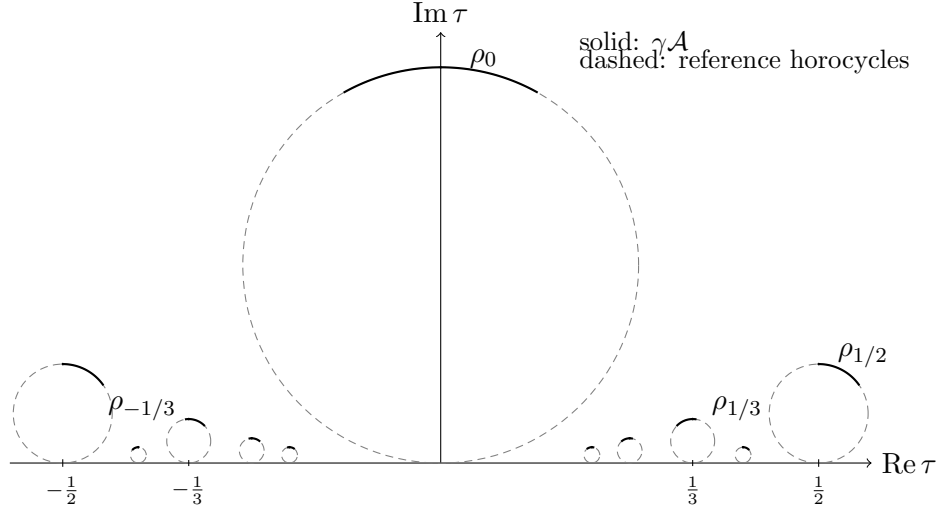

Figures~\ref{fig:real-Schwarz-arc} and \ref{fig:deformed-Ford} are
numerical plots of $h(\tau)=x$, obtained by Newton iteration.

\begin{lemma}\label{lem:trig-domination}
For every real $\theta$ and every positive integer $r$,
\[
 |\sin(r\theta)|\leq r|\sin\theta|.
\]
\end{lemma}

\begin{proof}
The assertion follows by induction from
\[
 \sin((r+1)\theta)
 =
 \sin(r\theta)\cos\theta+\cos(r\theta)\sin\theta.
\]
\end{proof}

\begin{proposition}
\label{prop:Watson-horizontal-sign}
For $x\notin\frac12\ZZ$,
\begin{equation}\label{eq:Watson-sign}
 \operatorname{sgn}\operatorname{Re}\mathcal G(x)
 =
 \operatorname{sgn}\sin(2\pi x),
 \quad
 \operatorname{sgn}\operatorname{Re}H(x)
 =
 -\operatorname{sgn}\sin(2\pi x).
\end{equation}
Consequently
\begin{equation}\label{eq:Watson-real-vanishing}
 \operatorname{Re}\mathcal G(x)=0
 \quad\Longleftrightarrow\quad
 \operatorname{Re}H(x)=0
 \quad\Longleftrightarrow\quad
 x\equiv0\ \text{or}\ \frac12\pmod1.
\end{equation}
\end{proposition}

\begin{proof}
Write
\[
 P_W(e^{2\pi ix})
 =
 \frac6\pi+\delta(x),
 \quad
 \delta(x)
 =
 \sum_{r\geq1}(-1)^{r-1}c_re^{2\pi irx}.
\]
Since the denominator is nonzero,
\[
 \operatorname{Re}\mathcal G(x)
 =
 \frac{\operatorname{Im}\delta(x)}
 {\left|6/\pi+\delta(x)\right|^2}.
\]
Watson's explicit majorant calculation in \cite[\S6]{Watson1952} gives positive
numbers $C_r$ with $|c_r|\leq C_r$, with
\[
 C_1=c_1>\frac{281630}{10^9}
\]
and
\[
 \sum_{r\geq1}rC_r<\frac{282726}{10^9}.
\]
Hence
\[
 \sum_{r\geq2}r|c_r|
 \leq
 \sum_{r\geq2}rC_r
 <
 \frac{1096}{10^9}
 <\frac{281630}{10^9}<c_1.
\]
Applying Lemma~\ref{lem:trig-domination} with $\theta=2\pi x$ gives
\[
 \left|
 \sum_{r\geq2}(-1)^{r-1}c_r\sin(2\pi rx)
 \right|
 <
 c_1|\sin(2\pi x)|.
\]
Thus $\operatorname{Im}\delta(x)$ has the sign of
$\sin(2\pi x)$ whenever the latter is nonzero. This proves the first
identity in \eqref{eq:Watson-sign}. The second follows from
\eqref{eq:G-reciprocal-H}, since
\[
 \operatorname{Re}\mathcal G(x)
 =
 -\frac{\operatorname{Re}H(x)}{|H(x)|^2}.
\]
The vanishing statement follows as well.
\end{proof}

\begin{proposition}
\label{prop:reference-horocycles}
Let
\[
 L_0=\left\{x+\frac{6i}{\pi}:x\in\RR\right\},
 \quad
 \gamma=\begin{pmatrix}a&b\\c&d\end{pmatrix},
 \quad c\neq0.
\]
Then $\gamma L_0$ is the Euclidean circle tangent to the real axis at $a/c$
with equation
\begin{equation}\label{eq:reference-circle}
 \left(u-\frac ac\right)^2
 +
 \left(v-\frac{\pi}{12c^2}\right)^2
 =
 \left(\frac{\pi}{12c^2}\right)^2.
\end{equation}
Its highest point is
\begin{equation}\label{eq:reference-circle-top}
 \frac ac+\frac{\pi i}{6c^2}.
\end{equation}
This is the principal term in the approximation to the zero near $a/c$
in \cite{ImamogluJermannToth2014}.
\end{proposition}

\begin{proof}
For $z=x+iY$ one has
\[
 \operatorname{Im}(\gamma z)
 =
 \frac{Y}{|cz+d|^2}.
\]
The image of the horizontal line $\operatorname{Im}z=Y$ is therefore the
circle tangent at $a/c$ with radius $1/(2c^2Y)$. Taking $Y=6/\pi$ gives
\eqref{eq:reference-circle}, and its highest point is twice the radius above
the real axis.
\end{proof}

\section{The CM case}

Masser determined the $\Qbar$-linear relations among periods and
quasi-periods of CM elliptic curves \cite{Masser1975,Masser1980}.

Assume first that $g_2,g_3\in\Qbar$ and that $E=\CC/\Lambda$ has complex multiplication. Write
\[
 \tau=\frac{\omega_2}{\omega_1}\in\HH.
\]
Then $\tau$ is imaginary quadratic. Masser's Theorem~III in \cite{Masser1975} gives dimension $4$ for the $\Qbar$-span of
\[
 1,\ \omega_1,\ \omega_2,\ \eta_1,\ \eta_2,\ 2\pi i,
\]
while Lemma~3.1 there gives the additional CM linear relation among the periods and quasi-periods. See also \cite[\S4.1]{Waldschmidt2006} for a modern summary. Besides $\omega_2=\tau\omega_1$, one may write the second relation in the form
\begin{equation}\label{eq:Masser-relation}
 C\eta_1-\tau\eta_2-\kappa\omega_2=0,
\end{equation}
where $C$ is algebraic and $\kappa\in\Qbar(\tau,g_2,g_3)$. Masser's analysis of the exceptional CM relation, in particular \cite[Lemma~3.1 and Theorem~III]{Masser1975} together with \cite{Masser1980}, shows that $\kappa=0$ precisely when $\tau$ is $\mathrm{SL}_2(\ZZ)$-equivalent to $i$ or to
\[
 \rho_3=e^{2\pi i/3}.
\]
See also \cite[\S4.1]{Waldschmidt2006}.

\begin{corollary}[Masser]\label{thm:Masser-CM-ratio}
Let $\Lambda=\ZZ\omega_1+\ZZ\omega_2$ be a CM lattice, and let $\eta_1,\eta_2$ be its quasi-periods. Then
\[
 [\eta_1:\eta_2]\in\PP^1(\Qbar)
 \quad\Longleftrightarrow\quad
 j(\Lambda)\in\{0,1728\}.
\]
Equivalently, the quasi-period ratio is algebraic exactly for the two $\mathrm{SL}_2(\ZZ)$-orbits of the elliptic fixed points $\rho_3$ and $i$.
\end{corollary}

\begin{proof}
The assertion is invariant under homothety, so a CM lattice may be rescaled so that $g_2,g_3\in\Qbar$. Masser's theorem \cite{Masser1975,Masser1980} then applies. If $\kappa=0$ in \eqref{eq:Masser-relation}, the relation gives an algebraic ratio between $\eta_1$ and $\eta_2$, and the characterization of the vanishing of $\kappa$ gives the two stated modular classes.

Conversely, suppose $r=\eta_1/\eta_2\in\Qbar$ and $\kappa\ne0$. Then
\[
 \eta_1-r\eta_2=0
\]
is a third $\Qbar$-linear relation among
$1,\omega_1,\omega_2,\eta_1,\eta_2,2\pi i$.
It is independent of the period relation $\omega_2-\tau\omega_1=0$ and of \eqref{eq:Masser-relation}: comparison of the $\eta_1,\eta_2$ coefficients would force a nonzero multiple of \eqref{eq:Masser-relation}, while comparison of the period coefficients would then force $\kappa=0$. This contradicts Masser's Theorem~III \cite{Masser1975}, according to which the six displayed quantities have $\Qbar$-span of dimension $4$, hence exactly two independent relations. Hence $\kappa=0$, and therefore $j\in\{0,1728\}$.
\end{proof}

\begin{corollary}\label{cor:CM-piE2}
If $\tau\in\HH$ is imaginary quadratic, then
\[
 \pi E_2(\tau)\in\Qbar
 \quad\Longleftrightarrow\quad
 j(\tau)\in\{0,1728\}.
\]
In particular, if $j(\tau)\notin\{0,1728\}$, then $\pi E_2(\tau)$ is transcendental.
\end{corollary}

\begin{proof}
For algebraic $\tau$, Proposition~\ref{prop:piE2} below identifies algebraicity of the quasi-period ratio with algebraicity of $\pi E_2(\tau)$. Theorem~\ref{thm:Masser-CM-ratio} gives the classification. At representatives of the two exceptional classes one has
\[
 \pi E_2(i)=3,
 \quad
 \pi E_2(\rho_3)=2\sqrt3,
\]
and algebraicity along their modular orbits also follows from the transformation formula \eqref{eq:piE2orbit} below.
\end{proof}

\section{Algebraic quasi-period ratios and the value of \texorpdfstring{$\pi E_2$}{pi E2}}

The corresponding non-CM assertion is conjectural.

\begin{conjecture}\label{conj:APQ}
If $\tau\in\Qbar\cap\HH$ is non-CM, then
\[
 [\eta_1:\eta_2]\notin\PP^1(\Qbar).
\]
\end{conjecture}

For algebraic $\tau$, the projective condition is equivalent to a statement about $\pi E_2(\tau)$.

\begin{proposition}\label{prop:piE2}
Let $\tau\in\Qbar\cap\HH$. Then
\[
 [\eta_1:\eta_2]\in\PP^1(\Qbar)
 \quad\Longleftrightarrow\quad
 \pi E_2(\tau)\in\Qbar.
\]
Consequently, for algebraic non-CM $\tau$, Conjecture~\ref{conj:APQ} is equivalent to the transcendence of $\pi E_2(\tau)$.
\end{proposition}

\begin{proof}
If $E_2(\tau)=0$, both assertions hold. Otherwise, from \eqref{eq:h},
\[
 \frac{\eta_2}{\eta_1}
 =\tau+\frac{6}{\pi iE_2(\tau)}.
\]
Since $\tau$ is algebraic, the quasi-period ratio is algebraic if and only if $\pi E_2(\tau)$ is algebraic. Projectivizing takes care of the cases in which one quasi-period vanishes.
\end{proof}

Proposition~\ref{prop:piE2} and Theorem~\ref{thm:Masser-CM-ratio} give the following global reformulation of Conjecture~\ref{conj:APQ}:
\begin{equation}\label{eq:sharp-piE2}
 \tau\in\Qbar\cap\HH
 \quad\Longrightarrow\quad
 \left(
 \pi E_2(\tau)\in\Qbar
 \ \Longleftrightarrow\ 
 j(\tau)\in\{0,1728\}
 \right).
\end{equation}
The CM case is Theorem~\ref{thm:Masser-CM-ratio}; the non-CM case
is Conjecture~\ref{conj:APQ}.

The condition is stable under the modular group. If
\[
 \gamma=\begin{pmatrix}a&b\\c&d\end{pmatrix}\in\mathrm{SL}_2(\ZZ),
\]
then
\begin{equation}\label{eq:piE2orbit}
 \pi E_2(\gamma\tau)
 =(c\tau+d)^2\pi E_2(\tau)+\frac{6c}{i}(c\tau+d).
\end{equation}
Hence, for algebraic $\tau$,
\[
 \pi E_2(\tau)\in\Qbar
 \quad\Longleftrightarrow\quad
 \pi E_2(\gamma\tau)\in\Qbar.
\]
In particular, if a hypothetical algebraic point $\rho$ satisfied $E_2(\rho)=0$, then
\[
 \pi E_2(\gamma\rho)=\frac{6c}{i}(c\rho+d)\in\Qbar
\]
for every $\gamma\in\mathrm{SL}_2(\ZZ)$.

If $[\eta_1:\eta_2]\in\PP^1(\QQ)$, Propositions
\ref{prop:rational-zero} and \ref{prop:jtrans} give
$j(\tau)\notin\Qbar$. Hence $g_2,g_3$ cannot both be algebraic in
any homothetic Weierstrass realization.

\section{Conditional period implications}
\label{sec:conditional-periods}

Let $E/K$ be an elliptic curve over an algebraically closed subfield of
$\CC$. In the Weierstrass basis $dx/y,x\,dx/y$, its period matrix is
\[
 \Pi_E=
 \begin{pmatrix}
  \omega_1&-\eta_1\\
  \omega_2&-\eta_2
 \end{pmatrix},
 \qquad \det\Pi_E=2\pi i.
\]
For $K_E=\aclC{\QQ(j(E))}$, let $\mathcal P(E/K_E)$ be the field
generated by $K_E$ and these periods. Fonseca's formulation of the
Grothendieck--Andr\'e period conjecture predicts
\begin{equation}\label{eq:AGPC}
 \trdeg_{\QQ}\mathcal P(E/K_E)\ge\dim\operatorname{MT}(E),
\end{equation}
with right-hand side $4$ for a non-CM elliptic curve
\cite[Conjecture~0.7 and Definition~13.1]{Fonseca2017}.

\begin{proposition}[Conditional consequence of Grothendieck--Andr\'e]
\label{thm:AG-implies}
Assume \eqref{eq:AGPC}. Then an algebraic non-CM
$\tau=\omega_2/\omega_1$ cannot have
$[\eta_1:\eta_2]\in\PP^1(\Qbar)$. Consequently every zero of $E_2$
is transcendental.
\end{proposition}

\begin{proof}
Since $\tau$ is algebraic and non-CM, Schneider's theorem makes
$j(\tau)$ transcendental.  Thus
$K_E=\aclC{\Qbar(j(\tau))}$ has transcendence degree one over
$\QQ$.  The algebraic period ratio gives
\[
 \omega_2=\tau\omega_1,
\]
so the two periods contribute at most one further generator over
$K_E$.  Likewise, the projective relation
$[\eta_1:\eta_2]\in\PP^1(\Qbar)$ leaves at most one quasi-period
generator: if $\eta_1\ne0$, write $\eta_2=r\eta_1$ with
$r\in\Qbar$, while if $\eta_1=0$ retain only $\eta_2$.  Consequently
the full period field is contained in a field of the form
$K_E(\omega_1,\eta)$, and hence
\[
 \trdeg_{\QQ}\mathcal P(E/K_E)\le1+2=3,
\]
contrary to the lower bound $4$ in \eqref{eq:AGPC}.

For the final assertion, suppose that an $E_2$-zero $\rho$ were
algebraic.  Proposition~\ref{thm:noCM} shows that it is non-CM, and
the normalized lattice $\ZZ+\ZZ\rho$ has $\eta_1=0$ by
\eqref{eq:eta-E2-normalized}.  It therefore satisfies the hypotheses
just ruled out.  Hence every $E_2$-zero is transcendental under
\eqref{eq:AGPC}.
\end{proof}

A different sufficient statement is the five-number conjecture,
recorded as Conjecture~1.11 in \cite{NesterenkoPhilippon2001}: for
non-CM $\tau$,
\begin{equation}\label{eq:five}
 \trdeg_{\QQ}\QQ\bigl(\tau,e^{2\pi i\tau},E_2(\tau),
 E_4(\tau),E_6(\tau)\bigr)\ge4.
\end{equation}
If an $E_2$-zero $\rho$ were algebraic, the field in
\eqref{eq:five} would be generated, up to algebraic elements, by the
three numbers $e^{2\pi i\rho},E_4(\rho),E_6(\rho)$. Nesterenko proves
that these three are algebraically independent, so its transcendence
degree would be exactly $3$, a contradiction.

For a zero $\rho$, the normalized lattice $\ZZ+\ZZ\rho$ has
\[
 \eta_1=0,\qquad \eta_2=-2\pi i,
\]
and therefore
\[
 \QQ\bigl(e^{2\pi i\rho},j(\rho),1,\rho,\eta_1,\eta_2\bigr)
 =\QQ\bigl(e^{2\pi i\rho},j(\rho),\rho,2\pi i\bigr),
\]
which has transcendence degree at most $4$. Thus the bound $5$ in
the analytic-lattice formulation of \cite[p.~196]{Diaz2000} cannot
hold for every normalization. The motivic period conjecture includes
the field of definition and the de Rham realization. We formulate the
following weaker bound for analytic lattices.

\begin{conjecture}
\label{conj:corrected-analytic-lattice}
Let $\tau\in\HH$ be non-CM, and let
$\Lambda=\ZZ\omega_1+\ZZ\omega_2$ be any lattice with
$\tau=\omega_2/\omega_1$ and quasi-periods $\eta_1,\eta_2$. Then
\begin{equation}\label{eq:corrected-analytic-lattice}
 \trdeg_{\QQ}\QQ\bigl(
 \tau,e^{2\pi i\tau},j(\tau),
 \omega_1,\omega_2,\eta_1,\eta_2
 \bigr)\ge4.
\end{equation}
\end{conjecture}

At an $E_2$-zero, the field in \eqref{eq:corrected-analytic-lattice}
for $\ZZ+\ZZ\rho$ is
\[
 \QQ\bigl(\rho,e^{2\pi i\rho},j(\rho),2\pi i\bigr),
\]
which has at most four generators. If $\rho$ were algebraic, it
would have transcendence degree at most three, contrary to the
conjecture.

\section{Waldschmidt's theorem and the rank-two obstruction}
\label{sec:waldschmidt-obstruction}

Waldschmidt's theorem \cite[Theorem~5.1]{Waldschmidt1983} concerns
analytic subgroups over a field of transcendence degree one.
For an analytic homomorphism $\varphi:\CC^n\to G(\CC)$ and a
subgroup $Y$ of rank $m$, let $d$ be the dimension of the Zariski
closure of the image. Set $p(G)=1$ for a linear group and $p(G)=2$
otherwise. The theorem requires
\begin{equation}\label{eq:Wald1983-numerical}
 md\ge2n(m+dp).
\end{equation}
Under this hypothesis and the other assumptions of the theorem, one
obtains a proper algebraic subgroup satisfying a stronger inequality.
For the homomorphism below, \eqref{eq:Wald1983-numerical} fails.

\begin{proposition}\label{prop:Waldschmidt-obstruction}
Assume that $\tau=\omega_2/\omega_1\in\Qbar\cap\HH$ is non-CM and that
\[
 r=\frac{\eta_2}{\eta_1}\in\Qbar.
\]
Let $K=\aclC{\Qbar(j(\tau))}$ and choose a Weierstrass model $E/K$. After a common homothety, choose a period basis $\omega_1,\omega_2$ of $E$ with $\omega_2/\omega_1=\tau$, and let $\eta_1,\eta_2$ be the corresponding quasi-periods. Define
\[
 \Phi:\CC\longrightarrow \mathbb G_a^2\times E(\CC),\quad
 z\longmapsto
 \left(\frac{z}{\omega_1},\frac{\eta_1 z}{2\pi i},\exp_E(z)\right).
\]
Then, for $Y=\ZZ\omega_1+\ZZ\omega_2$,
\[
 \Phi(Y)\subset (\mathbb G_a^2\times E)(K).
\]
Moreover the Zariski closure of $\Phi(\CC)$ over $K$ is the two-dimensional subgroup
\[
 H=\left\{(x,y,P):y=\frac{x}{\tau-r}\right\}\simeq \mathbb G_a\times E.
\]
Consequently, in Waldschmidt's notation, the natural data are
\[
 n=1,\quad m=2,\quad d=2,\quad p(H)=2,
\]
and the numerical hypothesis \eqref{eq:Wald1983-numerical} becomes $4\ge12$, which fails.
\end{proposition}

\begin{proof}
Legendre's relation gives
\[
 \omega_1\eta_1(\tau-r)=2\pi i,
\]
so $\tau-r\ne0$ and
\[
 \frac{\omega_1\eta_1}{2\pi i}=\frac1{\tau-r}\in\Qbar.
\]
For $z=m\omega_1+n\omega_2\in Y$ one therefore has
\[
 \Phi(z)=
 \left(m+n\tau,\frac{m+n\tau}{\tau-r},0_E\right),
\]
which is $K$-rational because $\tau,r\in\Qbar\subset K$.  This proves
$\Phi(Y)\subset(\mathbb G_a^2\times E)(K)$. The restriction
$\Phi|_Y$ is injective: $m+n\tau=0$ implies $m=n=0$ because
$\tau\notin\RR$. Thus $\Phi(Y)$ has rank two.

For arbitrary $z\in\CC$, the first two coordinates of $\Phi(z)$
satisfy the $K$-rational linear equation
\[
 y=\frac{x}{\tau-r}.
\]
Hence the Zariski closure $L$ of $\Phi(\CC)$ is an algebraic subgroup
of the displayed group $H\simeq\mathbb G_a\times E$.  Both
projections of $L$ are surjective.  Indeed, the first coordinate
$z/\omega_1$ ranges through all of $\CC$, while the elliptic
exponential $\exp_E:\CC\to E(\CC)$ is surjective.  The group $L$ is
connected: the composite of the connected group $\CC$ with the finite
component group $L/L^\circ$ is trivial, so $\Phi(\CC)\subset L^\circ$
and Zariski density forces $L=L^\circ$.  If $L$ were a proper
algebraic subgroup of $H$, it would therefore have dimension
one.  Surjectivity of both projections would then make $L$ isogenous
both to $\mathbb G_a$ and to $E$.  This is impossible: in
characteristic zero there is no nonzero algebraic homomorphism between
the additive group and an elliptic curve.  Therefore $L=H$ and
$d=\dim H=2$.

The source of $\Phi$ is one-dimensional, so $n=1$, and the preceding
injectivity gives $m=\operatorname{rank}\Phi(Y)=2$.  Since $H$ has an
elliptic factor it is not linear, and Waldschmidt's parameter is
$p(H)=2$.  Substitution into \eqref{eq:Wald1983-numerical} gives
\[
 md=2\cdot2=4,
 \qquad
 2n(m+dp)=2(2+2\cdot2)=12,
\]
so the required numerical hypothesis fails.
\end{proof}

For conjectural inequalities over finitely generated fields, see
\cite[\S6]{Waldschmidt1983}. Period conjectures that include the field
of definition are discussed in \cite{Andre2004,Bertolin2002}.

\section*{Acknowledgments} This work was partially supported by an NSERC Discovery Grant.

The author thanks Michel Waldschmidt for bringing his paper on analytic subgroups to his attention and for pointing out other references on the subject.

Figure 1 and Figure 2 were generated by ChatGPT.

\end{document}